\documentclass[10pt]{article}
\usepackage{amsfonts}
\usepackage{amsmath,amsthm,amssymb,mathrsfs,textcomp}
\usepackage{color}
\usepackage[hmargin=3cm,vmargin=3.5cm]{geometry}
\usepackage{graphicx}
\usepackage{cancel}
\usepackage{cite}
\usepackage{bm}

\newcommand{\ds}{\displaystyle}

\newcommand{\xb}{{\bf{x}}}

\newcommand{\R}{\mathbb{R}}

\newcommand{\la}{{\lambda}}

\newcommand{\vr}{{\varrho}}

\newcommand{\di}{{\rm div\, }}
\newcommand{\g}{{\nabla}}
\newcommand{\bu}{\mathbf u}

\theoremstyle{plain}
\newtheorem{theorem}{Theorem}[section]
\newtheorem{lemma}[theorem]{Lemma}

\newtheorem{remark}{Remark}[section]
\numberwithin{equation}{section} \numberwithin{theorem}{section}
\numberwithin{remark}{section} 
\input{tcilatex}
\begin{document}

\title{A Linearized Viscous, Compressible Flow-Plate Interaction \\ with Non-dissipative Coupling}
\author{George Avalos,\thanks{University of Nebraska-Lincoln, gavalos@math.unl.edu} \and Pelin G. Geredeli,\thanks{University of Nebraska-Lincoln, pguvengeredeli2@unl.edu; Hacettepe University, Ankara, Turkey, pguven@hacettepe.edu.tr} \and Justin T. Webster,\thanks{University of Maryland, Baltimore County, Maryland, websterj@umbc.edu} }
\maketitle

\begin{abstract}
We address semigroup well-posedness for a
linear, compressible viscous fluid interacting at its boundary with an elastic plate. We derive the model by linearizing the compressible Navier-Stokes equations about an arbitrary flow state,  so the fluid PDE includes an ambient flow profile $%
\mathbf{U}$. In contrast to model in \cite{AGW1}, we track the effect of this term at the flow-structure interface, yielding a velocity matching condition involving the material derivative of the structure; this destroys the dissipative nature of the coupling of the dynamics. We
adopt here a Lumer-Phillips approach, with a view of associating fluid-structure solutions with a $C_{0}$-semigroup $\left\{ e^{%
\mathcal{A}t}\right\} _{t\geq 0}$ on a chosen finite energy space of
 data. Given this approach, the  challenge becomes establishing the maximal dissipativity of an
operator $\mathcal{A}$, yielding the flow-structure dynamics.


\end{abstract}

\maketitle

\section{Introduction}

In mathematical studies of fluid-structure interactions arising in application the effect of viscosity can be important. Indeed, viscous fluids  introduce energy dissipation into the system, and produce non-trivial frictional effects in the interaction between fluid and solid. Interactive dynamics with viscous fluids are of paramount concern in the design and control of many physical systems \cite{dowell1}, e.g., aircraft, buildings and bridges, gas pipelines, engines, as well as other applications such as blood flow in an artery (for instance see \cite{canic} and references therein). In such applications, the density of the flow may change along a streamline, and \emph{compressibility}---the volume change  per unit pressure change---becomes non-negligible.

In many scenarios, mathematical solution techniques and analytical frameworks can be greatly simplified by assuming the flow is inviscid---viscosity free. One  of the principal motivating applications here is the field of {\em aeroelasticity}: elastic structures interacting with surrounding fluid flows, as with airfoils and paneling of aircraft.  Typically, the compressible gas is assumed to be inviscid and the flow to be
irrotational ({\em inviscid potential flow}) \cite{dowell1}. These assumptions reduce the flow dynamics  to a perturbed wave equation \cite{dcds,
dowell1,webster}. Yet there are situations
where viscous effects can simply not be neglected, e.g., for flows with Mach numbers\footnote{The Mach number $M$ is the ration of the flow velocity to the local speed of sound.}$ 0.3<M<3$ \cite{dowell1,fcomp}. Moreover, in certain regimes or configurations viscous effects are paramount, e.g., the low-speed flapping flag problem \cite{flag1,flag3} or the transonic flow regime $(0.8<M<1.2)$ \cite{dowell1}.

For incompressible flows, the analysis typically
involves two unknowns: velocity and pressure (e.g., \cite{visc2,canic,cr-full-karman,ChuRyz2011}). One must solve both conservation of mass and linear
momentum equations, with the fluid density constant. On the other hand, for
compressible flow, density and pressure vary in the dynamics.
Consequently, solutions in the compressible case require  
an equation of state for the fluid and a conservation
of energy statement. 
Here, we will assume the pressure depends linearly on the density. (For more discussion on this point, see Section \ref{modeling}.)

In this paper we begin with a laminar, unperturbed  {\em flow} of compressible fluid, and study perturbations about this given flow state. Such perturbations will be induced by a coupling to a non-stationary elastic dynamics imbedded in the fluid's boundary. We restrict our attention to such flow-structure interactions where the fluid exists in a 3-D spatial domain, bounded by a 2-D Lipschitz domain. A flat portion of the boundary is the equilibrium state of an elastic plate---with dynamics dictated by a fourth order plate equation that neglects the effects of rotational inertia. The flow and structure are  strongly coupled at the fluid-structure interface, with the plate dynamics affecting the flow through a normal component boundary condition, and the flow dynamics providing the dynamic  distributed stresses  across the surface of the plate. Taking the flow dynamics to be given by a linearization of compressible Navier-Stokes about a {\em non-zero} flow state $\mathbf U$ will turn out to have important repercussions in the interior (flow) dynamics, as well as at the flow-structure interface.

This problem is motivated by the main lines of the recent work of I. Chueshov et al. \cite{Chu2013-comp,cr-full-karman,ChuRyz2012-pois,ChuRyz2011} on various types of 3-D/2-D fluid-structure interactions (as considered here). In these papers, the authors address various combinations of viscous compressible and incompressible fluid dynamics in 3-D domains (sometimes unbounded but tubular \cite{ChuRyz2012-pois}) linearized about the steady flow state $\mathbf U\equiv \mathbf 0$, and coupled to different types of elastic dynamics at the interface (2-D in-plane elasticity \cite{berlin11}, von Karman \cite{ChuRyz2011}, or even full von Karman \cite{cr-full-karman}). The Galerkin approach is utilized in constructing solutions, and dynamical systems techniques are used to obtain long-time behavior results for these systems. As a primary motivating reference, in \cite{Chu2013-comp}, Chueshov considered the dynamics of a nonlinear plate, located on a flat portion of the boundary of a 3-D cavity, as it interacts with a compressible, isothermal
 fluid filling the cavity. There, the author addresses a natural first case of interest: $\mathbf U\equiv 0$, i.e., linearization about the trivial flow steady state. He shows both well-posedness and the existence of global attractors. In a personal correspondence with the authors of the present paper, Chueshov suggested that his method---based on a Galerkin approach with a priori energy estimates---would not accommodate the case of interest in aeroelasticity: linearization about $\mathbf U \neq \mathbf 0$. Thus, he suggested the problem at hand as a problem of interest, and intimated that a Lumer-Philips (semigroup) approach might yield well-posedness. 
 
In the present authors' previous work \cite{AGW1}, the suggested model was analyzed---that of \cite{Chu2013-comp}---{\em with  additional interior terms} associated to the $\mathbf U \neq 0$. For this non-dissipative flow-structure model, a pure velocity matching condition was imposed at the interface. {\em This type of coupling does not take into account flow effects at the interface with the plate} arising from the presence of  $\mathbf U \neq 0$. However, establishing well-posedness for the equations with a dissipative coupling already presented non-trivial technical challenges (discussed in detail in Section \ref{techreview}).

 Viewing the work in \cite{AGW1} as a requisite preliminary step, the  work at hand is the natural sequel  in justifying and mathematically accommodating the non-dissipative coupling. The present paper carefully derives the fluid-structure interface conditions, which {\em are necessarily non-dissipative} due to $\mathbf U\neq \mathbf 0$. Indeed, as with the fluid-structure models given in \cite{bolotin,webster,supersonic,btz}, we do not have a ``pure" velocity matching condition in the present work. Rather, we have an interface/coupling condition written also in terms of the material derivative\footnote{This condition arises via impermeability of the interface. See the calculation before \eqref{beforecalc}, or \cite[pp.172--174]{dowell1}.} of the structure, $(\partial_t+\mathbf U\cdot \nabla)w$. 
We critically use the techniques developed in \cite{AGW1} here, and the principal challenge is overcoming the addition of the non-dissipative coupling on the lower-dimensional interface. 
 
\subsection{Notation}
For the remainder of the text we write $\mathbf{x}$ for $(x_{1},x_{2},x_{3})%
\in \mathbb{R}_{+}^{3}$ or $(x_{1},x_{2})\in \Omega \subset \mathbb{R}%
_{\{(x_{1},x_{2})\}}^{2}$, as dictated by context. For a given domain $D$,
its associated $L^{2}(D)$ will be denoted as $||\cdot ||_D$ (or simply $%
||\cdot||$ when the context is clear). The symbols $\mathbf{n}$ and $\bf \tau$ will be used to denote, respectively, the unit external
normal and tangent vectors to $\mathcal{O}$. Inner products in $L^{2}(%
\mathcal{O})$ or $\mathbf{L}^{2}(\mathcal{O})$ are written $(\cdot ,\cdot)_{%
\mathcal{O}}$ (or simply $(\cdot ,\cdot)$ when the context is clear), while
inner products $L^{2}(\partial \mathcal{O})$ are written $\langle \cdot
,\cdot \rangle$. We will also denote pertinent duality pairings as $%
\left\langle \cdot ,\cdot \right\rangle _{X\times X^{\prime }}$, for a given
Hilbert space $X$. The space $H^{s}(D)$ will denote the Sobolev space of
order $s$, defined on a domain $D$, and $H_{0}^{s}(D)$ denotes the closure
of $C_{0}^{\infty }(D)$ in the $H^{s}(D)$-norm $\Vert \cdot \Vert
_{H^{s}(D)} $ or $\Vert \cdot \Vert _{s,D}$. We make use of the standard
notation for the boundary trace of functions defined on $\mathcal{O}$, which
are sufficently smooth: i.e., for a scalar function $\phi \in H^{s}(\mathcal{%
O})$, $\frac{1}{2}<s<\frac{3}{2}$, $\gamma (\phi )=\phi \big|_{\partial 
\mathcal{O}},$ a well-defined and surjective mapping on this range of $s$,
owing to the Sobolev Trace Theorem on Lipschitz domains (see e.g., \cite%
{necas}, or Theorem 3.38 of \cite{Mc}).

\section{Flow and Interface Modeling}\label{modeling}
Let $\mathcal{O}\subset \mathbb{R}^{3}$ be a \emph{bounded} and \emph{convex 
}fluid domain (and so has Lipschitz boundary $\partial \mathcal{O}$; see
e.g., Corollary 1.2.2.3 of \cite{grisvard}). The boundary decomposes into
two pieces $\overline{S}$ and $\overline{\Omega }$ where $\partial \mathcal{O%
}=\overline{S}\cup \overline{\Omega }$, with $S\cap \Omega =\emptyset $. We
consider $S$ to be the solid boundary, with no interactive dynamics, and $%
\Omega $ to be the equilibrium position of the elastic domain, upon which
the interactive dynamics take place. We also assume that: (i) the active
component $\Omega \subset \mathbb{R}^{2}$ is flat, with $C^2$
boundary, and embedded in the $x_{1}$--$x_{2}$ plane; (ii) the inactive
component $S$ lies below the $x_{1}$--$x_{2}$ plane. This is to say, 
\begin{align}
\Omega \subset & ~\{\mathbf{x}=(x_{1},x_{2},0)\} \\
S\subset & ~\{\mathbf{x}=(x_{1},x_{2},x_{3})~:~x_{3}\leq 0\}.
\end{align}%
We denote the unit outward normal vector to $
\partial \mathcal{O}$ by $\mathbf{n}(\mathbf{x})$ where $\left. \mathbf{n}\right\vert _{\Omega
}=[0,0,1],$ as in Figure 1 and the unit outward normal vector to $\partial \Omega$ by $\mathbf \nu(\mathbf x)$. 

\begin{figure}[tbp]
\begin{center}
\includegraphics[scale=.8]{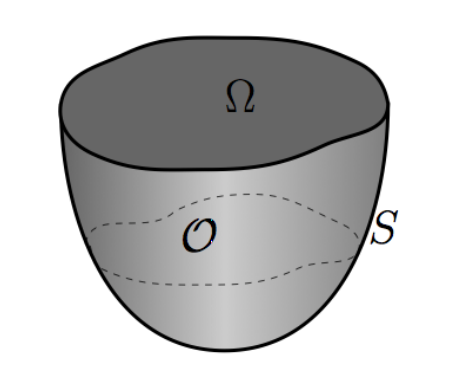}
\end{center}
\caption{The Fluid-Structure Geometry}
\end{figure}

Suppose  the domain $\mathcal O$ is filled with fluid whose governing dynamics
are the {\em compressible} Navier--Stokes
system \cite{fcomp,chorin-marsden} (see also \cite{berlin11}). We then  linearize  this model
with respect to a reference state $\{\vr_*;\mathbf U;p_*\}$, and suppose
that the unperturbed flow  $\mathbf{U}:\mathcal{O}\rightarrow \mathbb{R}^{3}$ is given by:

\begin{equation}
\mathbf{U}%
(x_{1},x_{2},x_{3})=[ U_{1}(x_{1},x_{2},x_{3}),U_{2}(x_{1},x_{2},x_{3}),U_{3}(x_{1},x_{2},x_{3})],
\label{flowfield}
\end{equation}%
and represents a mild (time-independent) ambient fluid flow. The quantities  $\vr_*,p_*$ are taken to be constant in time.
Then, for small perturbations $\{\rho; \bu; p\}$ of this ambient state, we write
\begin{align*}
\tilde \rho =~ \rho_*+\rho,~~~~~
\tilde p =~ p_*+p,~~~~~
\tilde{\bu} =~ \mathbf U+\bu.
\end{align*}
At this point, {\em we assume 
the pressure is a linear function of the density}. This assumption can be arrived at in two ways, which we now briefly describe\footnote{The discussion of pressure-density relations below was informed by personal correspondence with Earl Dowell \cite{dowellcorr}.}: 

For {\em isentropic} flows, the relationship between pressure and density is (see, e.g., \cite{fcomp,valli,dV})
\begin{equation}\label{baro} \tilde p = C[\tilde \rho]^{\gamma},\end{equation}
where $C$ is a constant evaluated for the pressure and density in the far field and $\gamma >0$---for air, $\gamma=1.4$. (We note that  {isentropic flows} \cite[pp.169--200]{dowell1} are barotropic---i.e., pressure depending only on density.)  Equation \eqref{baro} can be linearized through the above perturbation convention, taking $p_*$ and $\rho_*$ to be the far field pressure and density, yielding the linear relation
\begin{equation}\label{isen}{p}= C\big(\gamma,{\rho_*}, {p_{*}}\big){\rho}.\end{equation}
On the other hand, if one considers {\em isothermal} flow, the ideal gas law reads:
$\tilde p = \tilde \rho R T$, where $T$ is the temperature and $R$ is a fluid-dependent constant.
This equation also presents a linear relation between pressure and density, if $T$ is a constant.
Isothermal flows are used in low speed situations, i.e., with velocities much less than the speed of sound. Isentropic flow is used for compressible flows with small viscosity, and the ideal gas law is used for compressible viscous flows. We do note that $T$ is typically taken as an unknown for compressible viscous flows; if this consideration is made, an energy balance equation is required.

\begin{remark} Note that our assumption of linear pressure-density dependence could possibly be weakened in the analysis. In this case, we would take the flow to be isentropic as in \cite{valli,dV}, though, in these references the problem is entirely {\em stationary}. With this assumption, the flow coupling is inherently nonlinear, introducing additional complexity into the analysis \cite{fcomp}. We view the analyses in \cite{AGW1} and here as first steps, noting that the linear flow problem here is already of great technical complexity and mathematical challenge. We note the more recent references \cite{visc1,visc3}, each of which focuses on the issue of weak solvability for certain nonlinear, compressible, viscous fluid-structure interactions.
\end{remark} 

For mathematical simplicity we now take $p=\rho$ and assume $p_*=\vr_*=1$. If we generalize the forcing functions, then we obtain the physical {\em perturbation equations}:
\begin{subequations}\label{bar-model1-U*}
\begin{align}
 &
   (\partial_t+\mathbf U \cdot \nabla)p+\di\, \bu+(\text{div}~ \mathbf U)p=f(\mathbf x)\quad {\rm in}~~ \mathcal O
   \times\R_+,\label{den-eq1U*}
   \\[2mm]
 &
   (\partial_t+\mathbf U\cdot \nabla)\bu -\nu\Delta \bu -(\nu+\la)\g \di\, \bu + \g p+ \nabla \mathbf U \cdot \bu+(\mathbf U \cdot \nabla \mathbf U)p = \mathbf F(\mathbf x)\quad {\rm in}~ \mathcal O
   \times\R_+,\label{flu-eq1U*}
\end{align}
\end{subequations}
where the dynamic viscosity of the fluid is given by $\nu>0$, and $\lambda\ge0$ is Lam\'{e}'s first parameter (both of which would  vanish in the case of inviscid fluid). Given the {Lam\'{e} Coefficients}, the 
\textit{stress tensor} $\sigma$ of the fluid is defined as 
\begin{equation*}
\sigma (\mathbf{\mu })=2\nu \epsilon (\mathbf{\mu })+\lambda \lbrack
I_{3}\cdot \epsilon (\mathbf{\mu })]I_{3},
\end{equation*}%
where the \textit{strain tensor }$\epsilon $ is given by 
\begin{equation*}
\epsilon _{ij}({\bf \mu})=\dfrac{1}{2}\left( \frac{\partial \mu _{j}}{\partial
x_{i}}+\frac{\partial \mu _{i}}{\partial x_{j}}\right) \text{, \ }1\leq
i,j\leq 3.
\end{equation*}%
 With this notation it is easy to see that 
\begin{equation*}
\text{div}~\sigma (\mathbf{\mu })=\nu \Delta \mathbf{\mu }+(\nu +\lambda
)\nabla \text{div}(\mathbf{\mu }).
\end{equation*}

While the linearized interior terms are by now tractable \cite{AGW1}, we must supply the fluid equation with the correct boundary conditions on $\partial \mathcal O$ that will necessarily involve the plate's deflections $w$ on $\Omega \subset \partial \mathcal O$. The full system (with structural equations) will be discussed in detail in the following section; here, we impose the so called {\em impermeability condition} on $\Omega$, namely, that no fluid passes through the elastic portion of the boundary during deflection \cite{bolotin,dowell1}. 

Let $S(a_1,a_2,a_3)=0$ describe the interface in Lagrangian coordinates in $\mathbb R^3$; also let
$\xb=[ x_1,x_2,x_3]$ be the  Eulerian position inside $\mathcal O$. Then, letting $w(x_1,x_2,t)$ represent the transverse ($x_3$) displacement of the plate on $\Omega$, we have that
 $$S\big(x_1,x_2,x_3-w(x_1,x_2;t)\big)\equiv \mathscr S(x_1,x_2,x_3;t)=0,$$ describes the time-evolution of the boundary. The impermeability condition
 requires that the material derivative ($\partial_t+\tilde {\mathbf u}\cdot \nabla_{\mathbf x}$) vanishes on the deflected surface \cite{bolotin,chorin-marsden,dowell1}:
$$\big(\partial_t \mathcal +\tilde{\mathbf u} \cdot \nabla_{\mathbf x}\big) \mathscr S=0, ~~~~~\tilde{\mathbf u}=\mathbf u+\mathbf U$$
Applying the chain rule, we obtain
\begin{equation}\label{condish} \nabla_{\mathbf x} S\cdot [ 0,0,-w_t]+\mathbf U\cdot\big(\nabla_{\mathbf x}S+[ -S_{x_3}w_{x_1}, -S_{x_3}w_{x_2}, 0]\big)= -\mathbf u \cdot \big(\nabla_{\mathbf x}S+[ -S_{x_3}w_{x_1}, -S_{x_3}w_{x_2}, 0]\big).\end{equation}
We identify $\nabla_{\mathbf x}S$ as the normal to the deflected surface; {\em assuming small deflections} and restricting to $(x_1,x_2) \in \Omega$, we can identify $\nabla_{\mathbf x}S\big|_{\Omega}$ with $\mathbf n\big|_{\Omega}=[ 0, 0, 1]$. Making use of \eqref{condish}, imposing that $\mathbf U \cdot \mathbf n =0$ on $\partial \mathcal O$ (see \eqref{V_0} and discussion), and discarding quadratic terms, this relation allows us  to write for $(x_1,x_2) \in \Omega$:
$$\mathbf n\cdot [ 0, 0, w_t ] +\mathbf U \cdot [ w_{x_1},w_{x_2},0] =\mathbf u \cdot \mathbf n.$$
This yields the desired flow boundary condition for the dynamics:
\begin{equation}\label{beforecalc}
\mathbf u \cdot \mathbf n\big|_{\Omega} = w_t+\mathbf U\cdot \nabla w.
\end{equation}

\section{Main Flow-Structure PDE Model}
Deleting non-critical lower order and the benign inhomogeneous terms in \eqref{den-eq1U*}-\eqref{flu-eq1U*} (see also \cite{AGW1}), we
obtain the essential \emph{perturbation equations} to be studied below: 
\begin{align}
& \left\{ 
\begin{array}{l}
p_{t}+\mathbf{U}\cdot \nabla p+\text{div}~\mathbf{u}=0~\text{ in }~\mathcal{O}%
\times (0,\infty ) \\[.05cm] 
\mathbf{u}_{t}+\mathbf{U}\cdot \nabla \mathbf{u}-\text{div}~\sigma (\mathbf{u}%
)+\eta \mathbf{u}+\nabla p=0~\text{ in }~\mathcal{O}\times (0,\infty ) \\[.05cm] 
[\sigma (\mathbf{u})\mathbf{n}-p\mathbf{n}]\cdot \boldsymbol{\tau }=0~\text{
on }~\partial \mathcal{O}\times (0,\infty ) \\[.05cm] 
\mathbf{u}\cdot \mathbf{n}=0~\text{ on }~S\times (0,\infty ) \\[.05cm] 
\mathbf{u}\cdot \mathbf{n}=w_{t}+\mathbf{U}\cdot\nabla w~\text{ on }~\Omega \times (0,\infty )%
\end{array}%
\right.  \label{system1} \\
& \left\{ 
\begin{array}{l}
w_{tt}+\Delta ^{2}w+\left[ 2\nu \partial _{x_{3}}(\mathbf{u})_{3}+\lambda 
\text{div}(\mathbf{u})-p\right] _{\Omega }=0~\text{ on }~\Omega \times
(0,\infty ) \\[.05cm]  
w=\nabla w\cdot \mathbf \nu=0~\text{ on }~\partial \Omega \times
(0,\infty )%
\end{array}%
\right.  \label{IM2} \\
& 
\begin{array}{c}
\left[ p(0),\mathbf{u}(0),w(0),w_{t}(0)\right] =\left[ p_{0},\mathbf{u}%
_{0},w_{0},w_{1}\right] .%
\end{array}
\label{IC_2}
\end{align}%
Here, $p(t):\mathbb{R}^{3}\rightarrow \mathbb{R}$ and $\mathbf{u}(t):\mathbb{%
R}^{3}\rightarrow \mathbb{R}^{3}$ (pointwise in time) are given as the (Eulerian)
pressure and the fluid velocity field, respectively. The quantity $\eta >0$
represents a drag force of the domain on the viscous fluid. The function $w(t) : \mathbb{R}^{2} \to \mathbb R$ gives the transverse deflection of points $\mathbf x \in \Omega$, evolving according to a {\em plate equation}. In addition, the
given quantity $\mathbf{\tau }$ in (\ref{system1}) is in the space $%
TH^{1/2}(\partial \mathcal{O)}$ of tangential vector fields of Sobolev index
1/2; that is,%
\begin{equation}
\mathbf{\tau }\in TH^{1/2}(\partial \mathcal{O)=}\{\mathbf{v}\in \mathbf{H}^{%
\frac{1}{2}}(\partial \mathcal{O})~:~\mathbf{v}\cdot \mathbf{n}=0~\text{ on }%
~\partial \mathcal{O}\}.\footnote{%
See e.g., p.846 of \cite{buffa2}.}  \label{TH}
\end{equation}

\begin{remark}For convenience, we will refer to the fluid boundary condition on $\Omega$
\begin{equation} \label{coupling} \mathbf u \cdot \mathbf n =w_t+\mathbf U \cdot \nabla w 
\end{equation} 
as the {\em coupling condition}. 
\end{remark}

The boundary conditions invoked in (\ref{system1}) are the so-called \emph{impermeability-slip} conditions \cite%
{bolotin,chorin-marsden}. As discussed above, the physical interpretation is that no fluid
passes through the boundary, and that there is no fluid stress in the tangential direction 
$\tau$.

\section{Previous Literature, Present Approach, and Challenges} \label{techreview}

 In this paper we consider a possibly viscous
\emph{compressible} fluid flow in 3-D, interacting with a 2-D elastic structure. The recent works \cite{visc1,visc3} deal with the analogous system to \eqref{system1}--\eqref{IC2} in the fully nonlinear (isentropic case), but focus only on existence of weak solutions, whereas \cite{canic,visc2} deal with viscous incompressible fluids. In this work we focus on Hadamard well-posedness of an appropriate initial boundary value problem for a linearized system. In fact, beginning with compressible Navier-Stokes and linearizing, one can obtain
several related fluid-plate cases which are important from an applied point of view. Those most relevant to the analysis here:
(i) {\em incompressible fluid}: \cite{ChuRyz2011,cr-full-karman,ChuRyz2012-pois,berlin11} as well as \cite{clark,george1,george2}; (ii)
 {\em compressible fluid}: \cite{supersonic, webster,springer}. 
 
 The above works are thematically united in that the elastic
structure is 2-D, and evolves on the boundary of the 3-D fluid domain. The surveys \cite{berlin11,dcds} provide a nice overview of the
modeling, well-posedness, and long-time behavior results for the family of
dynamics described above.
In any analysis, allowing compressibility yields additional variables,
and, as a result, well-posedness is not obtained straightforwardly \cite{visc1,visc3}.  The essential difficulty lies in showing the range
condition of a generator, since one has to address this additional density/pressure
component. Such a variable cannot be readily eliminated, and therefore
accounts for an elliptic equation to be solved. 
Our approach here (developed in \cite{AGW1})
is based on the application of a static well-posedness
result in \cite{dV} (also see \cite{LaxPhil}).

As previously mentioned, the work \cite{Chu2013-comp} is the primary
motivating reference for the current analysis. The techniques used in \cite{Chu2013-comp}  are consistent with
those in \cite{ChuRyz2011,cr-full-karman}, namely, a
Galerkin procedure is implemented, along with good a priori
estimates, to produce solutions. As with many fluid-structure
interactions, the critical issue in \cite{Chu2013-comp} is the appearance of
ill-defined interface traces. In the incompressible case, one can
recover negative Sobolev trace regularity for the pressure at the
interface via properties of the Stokes' operator \cite{visc2,canic,cr-full-karman,ChuRyz2011}. However, in the viscous
compressible case this is no longer true. 

The semigroup approach, used by the present authors in \cite{AGW1}, does not
require the use of approximate solutions. We overcome the
difficulty of trace regularity issues by exploiting cancellations at the
level of solutions with data in the generator. In this way we do not have to
work component-wise on the dynamic equations, though we must work carefully
(and component-wise) on the static problem associated with maximality of the
generator. This fact, along with the merely Lipschitz nature of the flow-structure geometry, are the main technical hurdles overcome in \cite{AGW1}. In addition, \cite{preprint} addresses the long-time decay properties of that model.

The primary technical hurdles associated with the analysis here are now described: \begin{itemize} \item The presence of the ambient flow $\mathbf{U} \neq \mathbf 0$ in the modeling (i.e., in linearizing Navier-Stokes) introduces the term $\mathbf{U}\cdot \nabla p$ into the pressure equation, which {\bf
does not} represent a bounded perturbation of some straightforward interior dynamics. 
\item Additionally, as seen in the dissipativity calculation in Section \ref{disscalc}, the term $\mathbf U\cdot \nabla w$ in the coupling condition \eqref{coupling} at the interface $\Omega\subset\partial \mathcal O$ destroys the dissipative nature of the dynamics. This term---in the normal component of the flow---cannot be straightforwardly treated as a perturbation the dynamics in \cite{AGW1}. To do so, one would need to have tight control of the dynamic boundary-to-interior mapping for the flow dynamics. This type of approach, for instance, is utilized in \cite{supersonic}, but there it requires a good structure of Neumann-lift maps for hyperbolic dynamics, and a viable dynamic trace regularity theory---neither of which are readily available here.  \end{itemize}

Thus, we must first alter the domain of the generator (conditions (A.i)--(A.iv) below) to accommodate the coupling conditions  \eqref{coupling}. Having made this choice, if we utilize the natural topology on $\mathcal H$ induced by \eqref{stand} (defined below), the associated terms that fall out of the dissipativity calculation ~$\text{Re}\left(\mathcal A\mathbf y , \mathbf y\right)_{\mathcal H, \text{stand.}}$  {\em involve traces of the fluid well above the energy level}. Thus, to induce dissipativity of the dynamics, we alter the inner-product structure on $\mathcal H$ in accordance with the ``adjustment" of $D(\mathcal A)$ (as was done in \cite{webster,graber}). (Formally, we change in the inner-product structure on $\mathcal H$ as motivated by the change of variable $w_t \mapsto D_t w= (\partial_t +\mathbf U\cdot \nabla)w$.)   

With the dynamics operator appropriately adjusted, and the state space topology accordingly changed, we proceed to estimate the resulting terms with an eye to obtain {\em shifted}-dissipativity---a bounded perturbation of the dynamics operator will be shown to be maximal dissipative. The requisite estimates follow from carefully looking at the $\mathbf U \cdot \nabla w$ term as if it were the boundary flux multiplier $h \cdot \nabla w$, used frequently in wave and plate dynamics (see, e.g., \cite[pp.574--579]{springer} and \cite{lagnese} and references therein) to obtain the equipartition of energy. We note that if one simply considers the multiplier $\mathbf U\cdot \nabla w$ in the equations, it will be necessary to control a term of the form $\ds (\mathbf U\cdot \mathbf \nu)\int_{\partial \Omega} |\Delta w|^2d\partial\Omega$; in practice, however, we have no control over the interaction with $\mathbf U\big|_{\Omega}$ and $\mathbf \nu$ on $\partial \Omega$. Thus our choice of inner-product on $\mathcal H$ must also be constructed so as to relegate such terms to ``lower order" status, so they can be absorbed via a {\em bounded perturbation}. 

\section{Functional Setup and the Main Result}

\label{results} We are primarily interested in Hadamard well-posedness of the
linearized coupled system given in (\ref{system1})--(\ref{IC_2}).
Specifically, we will ascertain well-posedness of the PDE model (\ref%
{system1})--(\ref{IC_2}) for arbitrary initial data in the natural space of
finite energy. To accomplish this, we will adopt a semigroup approach;
namely, we will pose and validate an explicit semigroup generator
representation for the fluid-structure dynamics (\ref{system1})--(\ref{IC_2}%
), yielding {\em strong} and {\em mild} solutions for the coupled system \cite{pazy}.

 For convenience, we now define the space 
\begin{equation}\label{V_0}
\mathbf{V}_{0}=\{\mathbf{v}\in \mathbf{H}^{1}(\mathcal{O})~:~\left. \mathbf{v%
}\right\vert _{\partial \mathcal{O}}\cdot \mathbf{n}=0~\text{ on }~\partial 
\mathcal{O}\};  
\end{equation}
With respect to the \textquotedblleft ambient flow\textquotedblright\ field $%
\mathbf{U}$ appearing in (\ref{system1}), we impose the assumptions: 
\begin{enumerate}
\item[(i)]
$\mathbf{U}\in \mathbf{V}_{0}\cap \mathbf{H}^{3}(\mathcal{O});$
\item[(ii)] $\mathbf{U}\big|_{\Omega} \in C^{2}(\overline{\Omega}).$
\end{enumerate}
\begin{remark}
For assumption (i), see the analogous---and actually slightly stronger---specifications made on
ambient fields on p.529 of \cite{dV} and pp.102--103 of \cite{valli}. Assumption (ii) was not present in the authors' recent work \cite{AGW1}, but is necessary to address the boundary term $\mathbf U \cdot \nabla w$ in the coupling condition.
\end{remark}

With respect to the coupled PDE system (\ref{system1})--(\ref{IC_2}), the
associated space of well-posedness will be 
\begin{equation}
\mathcal{H}\equiv [p,\bu, w_0,w_1] \in L^{2}(\mathcal{O})\times \mathbf{L}^{2}(\mathcal{O})\times
H_{0}^{2}(\Omega )\times L^{2}(\Omega). \label{H}
\end{equation}%

\noindent In what follows, we consider the linear operator $\mathcal{A}:D(\mathcal{A}%
)\subset \mathcal{H}\rightarrow \mathcal{H}$, which expresses the
 PDE system (\ref{system1})--(\ref{IC_2}) as the
abstract ODE: 
\begin{eqnarray}
\dfrac{d}{dt}%
\begin{bmatrix}
p \\ 
\mathbf{u} \\ 
w \\ 
w_{t}%
\end{bmatrix}
&=&\mathcal{A}%
\begin{bmatrix}
p \\ 
\mathbf{u} \\ 
w \\ 
w_{t}%
\end{bmatrix}%
;  \notag \\
\lbrack p(0),\mathbf{u}(0),w(0),w_{t}(0)] &=&[p_{0},\mathbf{u}%
_{0},w_{0},w_{1}] \in \mathcal{H}.  \label{ODE}
\end{eqnarray}

\noindent To wit, the action of $\mathcal A$ is given by
\begin{equation}
\mathcal{A}=\left[ 
\begin{array}{cccc}
-\mathbf{U}\mathbb{\cdot }\nabla (\cdot ) & -\text{div}(\cdot ) & 0 & 0 \\ 
-\mathbb{\nabla (\cdot )} & \text{div}\sigma (\cdot )-\eta I-\mathbf{U}%
\mathbb{\cdot \nabla (\cdot )} & 0 & 0 \\ 
0 & 0 & 0 & I \\ 
\left. \left[ \cdot \right] \right\vert _{\Omega } & -\left[ 2\nu \partial
_{x_{3}}(\cdot )_{3}+\lambda \text{div}(\cdot )\right] _{\Omega } & -\Delta
^{2} & 0%
\end{array}%
\right], \label{AAA}
\end{equation}

\noindent with the domain $D(\mathcal{A})$ given as 
\begin{equation*}
D(\mathcal{A})=\{(p_{0},\mathbf{u}_{0},w_{0},w_{1})\in L^{2}(\mathcal{O}%
)\times \mathbf{H}^{1}(\mathcal{O})\times H_{0}^{2}(\Omega )\times
H_{0}^{2}(\Omega )~:~(A.i)\text{--}(A.v)~~\text{hold below}\},
\end{equation*}%
where

\begin{enumerate}
\item[(A.i)] $\mathbf{U}\cdot \nabla p_{0}\in L^{2}(\mathcal{O})$;

\item[(A.ii)] $\text{div}~\sigma (\mathbf{u}_{0})-\nabla p_{0}\in L^{2}(%
\mathcal{O})$  (and so from this and an integration by parts, we also have 
\newline  
$\left[ \sigma (\mathbf{u}_{0})\mathbf{n}-p_{0}\mathbf{n}\right] _{\partial 
\mathcal{O}}\in \mathbf{H}^{-\frac{1}{2}}(\partial \mathcal{O})$).

\item[(A.iii)] $-\Delta ^{2}w_{0}-\left[ 2\nu \partial _{x_{3}}(\mathbf{u}%
_{0})_{3}+\lambda \text{div}(\mathbf{u}_{0})\right] _{\Omega }+\left.
p_{0}\right\vert _{\Omega }\in L^{2}(\Omega )$

\item[(A.iv)] $\left( \sigma (\mathbf{u}_{0})\mathbf{n}-p_{0}\mathbf{n}%
\right) \bot ~TH^{1/2}(\partial \mathcal{O})$. That is, 
\begin{equation*}
\left\langle \sigma (\mathbf{u}_{0})\mathbf{n}-p_{0}\mathbf{n},\mathbf{\tau }%
\right\rangle _{\mathbf{H}^{-\frac{1}{2}}(\partial \mathcal{O})\times 
\mathbf{H}^{\frac{1}{2}}(\partial \mathcal{O})}=0\text{ \ for all }\mathbf{%
\tau }\in TH^{1/2}(\partial \mathcal{O}).
\end{equation*}

\item[(A.v)] $\mathbf{u}_{0}=\mathbf{\mu }_{0}+\widetilde{\mathbf{\mu }}_{0}$, where $\mathbf{\mu }_{0}\in \mathbf{V}_{0}$ and $\widetilde{\mathbf{\mu }}%
_{0}\in \mathbf{H}^{1}(\mathcal{O})$ satisfies\footnote{%
The existence of an $\mathbf{H}^{1}(\mathcal{O})$-function $\widetilde{\mathbf{\mu }}_{0}$ with such a boundary trace on Lipschitz domain $\mathcal{%
O}$ is assured; see e.g., Theorem 3.33 of \cite{Mc}, or see also the proof
of Lemma \ref{staticwellp} below.}
\begin{equation*}
\left. \widetilde{\mathbf{\mu }}_{0}\right\vert _{\partial \mathcal{O}}=%
\begin{cases}
0 & ~\text{ on }~S \\ 
[w_{1}+\mathbf{U}\cdot \nabla w_0]\mathbf{n} & ~\text{ on}~\Omega,%
\end{cases}%
\end{equation*}%
\noindent (Note that $\left. \mathbf{\mu }_{0}\right\vert _{\partial \mathcal{O}%
}\in TH^{1/2}(\partial \mathcal{O})$ and since $w_1 +\mathbf U\cdot \nabla w_0 \in H_0^1(\Omega)$,  $\widetilde{\mathbf{\mu}}_0$ is well-defined).
\end{enumerate}

In the following theorem, we provide the semigroup well-posedness for $\mathcal{A%
}:D(\mathcal{A})\in \mathcal{H}\rightarrow \mathcal{H}$, the proof of which
is based on the well known Lumer-Phillips Theorem and associated bounded perturbation result \cite{pazy}.

\begin{theorem}
\label{wellp} The map $\left[p_{0},\mathbf{u}_{0};w_{0},w_{1}\right]
\rightarrow \left[ p(t),\mathbf{u}(t);w(t),w_{t}(t)\right]$ defines a
strongly continuous semigroup $\{e^{\mathcal{A}t}\}$ on the space $\mathcal{H%
}$, and hence the PDE system (\ref{system1})--(\ref{IC_2})---or equivalently, the initial value problem (\ref{ODE})---is well-posed. 

In particular:

\begin{enumerate}
\item[(i)]  If $\left[ p_{0},\mathbf{u}_{0},w_{0},w_{1}\right] \in D(%
\mathcal{A})$, then $\left[ p,\mathbf{u},w,w_{t}\right] \in C([0,\infty );D(%
\mathcal{A}))\cap C^{1}([0,\infty );\mathcal{H})$;

\item[(ii)] If $\left[ p_{0},\mathbf{u}_{0},w_{0},w_{1}\right] \in %
\mathcal{H}$, then $\left[ p,\mathbf{u},w,w_{t}\right] \in C([0,\infty );%
\mathcal{H})$. 
\[
\]
\end{enumerate}

In addition, the
semigroup enjoys the following estimate for some $\varepsilon>0$ sufficiently large: 
\begin{equation}
\big|\big|e^{\mathcal{A}t}\big|\big|_{\mathcal{L}(\mathcal{H})}\leq e^{Kt},~~\forall t>0,
\end{equation}
with $\ds K=\dfrac{1}{2}||\text{div}(\mathbf{U})||_{L^{\infty}(\mathcal O)}+\varepsilon$.
\end{theorem}

\begin{remark}
\label{weaker} Given the existence of a semigroup $\{e^{\mathcal{A}t}\}$ for
the fluid-structure generator $\mathcal{A}:D(\mathcal{A})\subset \mathcal{H}%
\rightarrow \mathcal{H}$: if initial data $[p_{0},\mathbf{u}%
_{0};w_{0},w_{1}]\in D(\mathcal{A})$, the corresponding solution $[p(t),%
\mathbf{u}(t);w(t),w_{t}(t)]\in C([0,\infty ),D(\mathcal{A}))$. In
particular, the solution satisfies the condition (A.iv) in the definition of
the generator. This means that one has the tangential boundary
condition 
\begin{equation*}
\lbrack \sigma (\mathbf{u}_{0})\mathbf{n}-p_{0}\mathbf{n}]\cdot \mathbf{\tau 
}=0\text{\ \ for all }\mathbf{\tau }\in TH^{1/2}(\partial \mathcal{O}),
\end{equation*}%
satisfied in the sense of distributions. That is to say, $\forall ~
\boldsymbol{\tau }\in TH^{1/2}(\partial \mathcal{O})$ and $\forall \phi \in 
\mathcal{D}(\partial \mathcal{O})$,%
\begin{equation}
\langle \sigma (\mathbf{u}_{0})\mathbf{n}-p_{0}\mathbf{n},\phi \boldsymbol{%
\tau }\rangle _{\partial \mathcal{O}}=0.  \label{weak2}
\end{equation}
\end{remark}

\section{The Proof of Theorem \protect\ref{wellp}}
\label{proof}

Our proof of well-posedness hinges on showing that the operator $\mathcal{A}:D(%
\mathcal{A})\subset \mathcal{H}\rightarrow \mathcal{H}$ generates a $C_{0}$%
-semigroup. As discussed in Section \ref{techreview}, the presence of a
generally nonzero ambient vector field $\mathbf{U}$ produces a
lack of dissipativity of the operator $\mathcal{A}$. Accordingly, we
introduce the following bounded perturbation $\widehat{\mathcal{A}}$ of the
generator $\mathcal{A}$: 
\begin{equation}
\widehat{\mathcal{A}}=\mathcal{A}-
\begin{bmatrix}
 \dfrac{\text{div}(\mathbf{U})}{2}+ \varepsilon I & 0 & 0 & 0 \\ 
0 & \dfrac{\text{div}(\mathbf{U})}{2}+ \varepsilon I & 0 & 0 \\ 
0 & 0 & \varepsilon I & 0 \\ 
0 & 0 & 0 & \varepsilon I%
\end{bmatrix}%
\text{, \ \ }D(\widehat{\mathcal{A}})=D(\mathcal{A}).  \label{pertA}
\end{equation}%
Therewith, the proof of Theorem \ref{wellp} is geared towards establishing
the maximal dissipativity of the linear operator $\widehat{\mathcal{A}}$;
subsequently, an application of the Lumer-Phillips Theorem will yield that $%
\widehat{\mathcal{A}}$ generates a $C_{0}$ semigroup of contraction on $%
\mathcal{H}$. In turn, applying the standard perturbation result \cite{kato} or \cite[Theorem 1.1, p.76]{pazy}, yields semigroup
generation for the original modeling fluid-structure operator $\mathcal{A}$
of (\ref{AAA}), via (\ref{pertA}).

\subsection{Inner Product and Induced Norm}\label{dissnorm}

 In what follows, we will require an equivalent modification of the natural
norm for finite energy space $\mathcal{H}$. To construct it, we will require three ingredients:

{(i)}
Define the \textquotedblleft Dirichlet map\textquotedblright, which
extends essential boundary data $\phi$ on $\Omega $ to a harmonic function in $%
\mathcal{O}$:
\begin{equation}
D\phi=f\Leftrightarrow \left\{ 
\begin{array}{l}
\Delta f=0\text{ \ in }\mathcal{O} \\ 
\left. f\right\vert _{\partial \mathcal{O}}=\phi _{ext}\text{ \ on }\partial 
\mathcal{O},
\end{array}%
\right.   \label{D}
\end{equation}%
where%
\begin{equation}
\phi _{ext}=\left\{ 
\begin{array}{l}
0\text{ \ on \ }S \\ 
\phi\text{ \ on }\Omega .%
\end{array}%
\right.   \label{ext}
\end{equation}%
Therewith, if $\phi \in H_{0}^{\frac{1}{2}+\epsilon }(\Omega )$, then $\phi
_{ext}\in H^{\frac{1}{2}+\epsilon }(\partial \mathcal{O})$ (see e.g.,
Theorem 3.3, p. 95, of \cite{Mc}). Subsequently, via Lax-Milgram, we deduce
that 
\begin{equation}
D\in \mathcal{L}\left( H_{0}^{\frac{1}{2}+\epsilon }(\Omega ),H^{1}(\mathcal{%
O})\right). \label{D1}
\end{equation}

{(ii)}  Let $\mathbf g(\xb)$ be a $C^2$-extension of $\mathbf \nu(\xb)$ (the unit normal to the boundary of $\Omega$). That is, 
\begin{equation}
\mathbf{g}\in C^{2}(\overline{\Omega })\text{, \ }\left. \mathbf{g}%
\right\vert _{\partial \Omega }={\mathbf \nu} (\mathbf{x}).  \label{g}
\end{equation}

{(iii)}  With $\mathbf g$ as above, we construct a function
\begin{equation}
\mathbf{h}_{\alpha}(\cdot )\equiv \left. \mathbf{U}\right\vert _{\Omega }(\cdot
)-\alpha \mathbf{g}(\cdot ),  \label{h}
\end{equation}
where $\alpha>0$ is a parameter that will eventually be taken sufficiently large.\\ \\
\vspace{0.5cm}
Now, with the ingredients (i) $D(\cdot)$, (ii) $\mathbf g(\cdot)$, (iii) $\mathbf {h}_{\alpha}(\cdot)$ in hand, and with unit vector $\mathbf{e}_3=[0,0,1]$, we  
topologize $\mathcal H$ inner-product in the following way: 
\begin{align}\label{innerproduct}
(\mathbf{y}_{1},\mathbf{y}_{2})_{\mathcal{H}}=&~(p_{1},p_{2})_{L^{2}(\mathcal{O
})} +(\mathbf{u}_{1}-\alpha D(\mathbf g\cdot \nabla w_1) \mathbf{e}_3,\mathbf{u}_{2}-\alpha D(\mathbf g\cdot \nabla w_2) \mathbf{e}_3)_{\mathbf{L}^{2}(\mathcal{O})} \\\nonumber
&+(\Delta
w_{1},\Delta w_{2})_{L^{2}(\Omega )}   +(v_{1}+\mathbf{h}_{\alpha}\cdot \nabla w_1,v_{2}+\mathbf{h}_{\alpha}\cdot \nabla w_2)_{L^{2}(\Omega )}
\end{align}%
for any $\mathbf{y}_{1}=[p_{1},\mathbf{u}_{1},w_{1},v_{1}] \in \mathcal{H}$ and $\mathbf{y}_{2}=[p_{2},\mathbf{u}_{2},w_{2},v_{2}] \in \mathcal{H}$. (Note that by \eqref{D}-\eqref{D1}, this inner product is well-defined.) This inner product induces a norm on $\mathcal H$, given by: for $\mathbf y = [p, \mathbf u, w, v]$
\begin{align}
||\mathbf y||^2_{\mathcal H}= ||p||^2_{L^2(\mathcal O)}+ ||\mathbf u-\alpha D(\mathbf g\cdot \nabla w) \mathbf{e}_3||^2_{\mathbf L^2(\mathcal O)}+||\Delta w||_{L^2(\Omega)}^2+||v+\mathbf{h}_{\alpha}\cdot \nabla w||_{L^2(\Omega)}^2.
\end{align}
We also note that the standard norm on $\mathcal H$ is given by 
\begin{align}\label{stand}
||\mathbf y||^2_{\mathcal H,\text{stand.}}= ||p||^2_{L^2(\mathcal O)}+ ||\mathbf u||^2_{\mathbf L^2(\mathcal O)}+||\Delta w||_{L^2(\Omega)}^2+||v||_{L^2(\Omega)}^2.
\end{align}
It is clear (by \eqref{D1}) that, since the appended terms in the modified inner product are lower order,  the norm induced by the inner product \eqref{innerproduct} on $\mathcal H$ is {\em equivalent} to $||\cdot||_{\mathcal H,\text{stand}}$, which was utilized in \cite{AGW1}.

\subsection{Dissipativity of $\widehat{\mathcal{A}}$}

\subsubsection{A Regularity Result for the Mapping $D$}

In the course of establishing the dissipativity of $\widehat{\mathcal{A}}:D(%
\mathcal{A})\rightarrow \mathcal{H}$, we will have to apply the Dirichlet
map $D$ of (\ref{D}) to boundary data rougher than $H_{0}^{\frac{1}{2}%
+\epsilon }(\Omega )$ (see \ref{D1}). Since the flow domain $\mathcal{O}$ is
Lipschitz, we cannot apply the known regularity results for second order
elliptic boundary value problems with rough boundary data in \cite[Remark 7.2, p. 188]{L-M}, 
which \textit{de facto} require $\mathcal{O}$ to have a $C^{2}$
boundary. Accordingly, we need the following lemma.

\begin{lemma}
\label{mark}The Dirichlet map $D$, as defined in (\ref{D}), is an element of $%
\mathcal{L}(H^{-\frac{1}{2}}(\Omega ),L^{2}(\mathcal{O})).$
\end{lemma}

\begin{proof}[Proof of Lemma \ref{mark}]{\bf Step 1. }We start by defining the Dirichlet Laplacian, 
$A_{D}:D(A_{D})\subset L^{2}(\mathcal{O})\rightarrow L^{2}(\mathcal{O})$ by 
\begin{equation}
A_{D}v=-\Delta v,\ D(A_{D})=H^{2}(\mathcal{O})\cap H_{0}^{1}(\mathcal{O}).
\label{A_D}
\end{equation}%
As defined, $A_{D}$ is a positive definite, self-adjoint operator. Moreover,
since the domain $\mathcal{O}$ is bounded and convex, then by \cite[Theorem 3.2.1.2,
p. 147]{grisvard}, we have that 
\[
A_{D}\text{ \textit{is an isomorphism from} }D(A_{D})\textit{ \text{onto} }%
L^{2}(\mathcal{O})\text{.}
\]%
In turn, by duality, we have that%
\begin{equation}
A_{D}\text{ \textit{is an isomorphism from} }L^{2}(\mathcal{O})\text{ 
\textit{onto }}\left[ D(A_{D})\right] ^{\prime }.  \label{iso}
\end{equation}

\medskip

{\bf Step 2.} Since $\mathcal{O}$ has Lipschitz boundary, then for given $%
v\in D(A_{D})$, its normal derivative $\left. \frac{\partial v}{\partial 
\mathbf{n}}\right\vert _{\Omega }$ is only assured to be in $L^{2}(\partial 
\mathcal{O})$; see \cite{necas}. However, from the relatively recent result
in \cite{buffa}, we have that for any $f\in H^{2}(\mathcal{O})$,%
\begin{equation}
\nabla _{\partial \mathcal{O}}f+\frac{\partial f}{\partial \mathbf{n}}%
\mathbf{n}\in \mathbf{H}^{\frac{1}{2}}(\partial \mathcal{O}),  \label{mem}
\end{equation}%
where $\left[ \nabla _{\partial \mathcal{O}}f\right] _{\partial \mathcal{O}}$
denotes the tangential gradient of $\left. f\right\vert _{\partial \mathcal{O%
}}$ (see Theorem 5 of \cite{buffa}).  Since $\left. \mathbf{n}\right\vert
_{\Omega }=[0,0,1]$, we then infer that, in particular, 
\begin{equation}
\left. \frac{\partial v}{\partial \mathbf{n}}\right\vert _{\Omega }\in H^{%
\frac{1}{2}}(\Omega )\text{, for every }v\in D(A_{D}).  \label{hid}
\end{equation}

{\bf Step 3.} Given boundary function $\phi \in L^{2}(\Omega )$, we denote
its extension by zero, as in (\ref{ext}), via%
\[
\phi _{ext}=\left\{ 
\begin{array}{c}
0\text{, \ on }S \\ 
\phi \text{, \ on }\Omega .%
\end{array}%
\right. 
\]%
Therewith, we define the linear functional $l_{\phi }$, by having for
any\thinspace $v\in D(A_{D})$, 
\begin{eqnarray}
l_{\phi }(v) &=&\left( \phi _{ext},\frac{\partial v}{\partial \mathbf{n}}%
\right) _{\partial \mathcal{O}}=\left( \phi ,\frac{\partial v}{\partial 
\mathbf{n}}\right) _{\Omega }  \nonumber \\
&=&\left\langle \phi ,\frac{\partial v}{\partial \mathbf{n}}\right\rangle
_{H^{-\frac{1}{2}}(\Omega )\times H^{\frac{1}{2}}(\Omega )},  \label{e1}
\end{eqnarray}%
where in the last equality, we are using (\ref{hid}) and the fact that $H^{%
\frac{1}{2}}(\Omega )=H_{0}^{\frac{1}{2}}(\Omega )$ (see e.g., of \cite[Theorem 3.40,
p. 105]{Mc}). An estimation of right hand side, via (\ref{hid}) and
the Closed Graph Theorem, then yields for given $\phi \in L^{2}(\Omega )$, 
\begin{equation}
\left\vert l_{\phi }(v)\right\vert \leq C\left\Vert \phi \right\Vert _{H^{-%
\frac{1}{2}}(\Omega )}\left\Vert v\right\Vert _{D(A_{D})}.  \label{e2}
\end{equation}%
A subsequent extension by continuity yields that for any $\phi \in H^{-\frac{%
1}{2}}(\Omega )$, that $l_{\phi }$, as given by 
\begin{equation}
l_{\phi }(v)=\left( \phi _{ext},\frac{\partial v}{\partial \mathbf{n}}%
\right) _{\partial \mathcal{O}}=\left\langle \phi ,\frac{\partial v}{%
\partial \mathbf{n}}\right\rangle _{H^{-\frac{1}{2}}(\Omega )\times H^{\frac{%
1}{2}}(\Omega )}\text{ \ }\forall \text{ }v\in D(A_{D})\text{,}  \label{e3}
\end{equation}%
is an element of $\left[ D(A_{D})\right] ^{\prime }$, with the estimate%
\begin{equation}
\left\Vert l_{\phi }\right\Vert _{\left[ D(A_{D})\right] ^{\prime }}\leq
C\left\Vert \phi \right\Vert _{H^{-\frac{1}{2}}(\Omega )}.  \label{e4}
\end{equation}

{\bf Step 4. }Via transposition, the problem of finding $f\in L^{2}(%
\mathcal{O})$ which solves the boundary value problem (\ref{D}) with
boundary data $\phi \in H^{-\frac{1}{2}}(\Omega )$, is the problem of
finding $f\in L^{2}(\mathcal{O})$ which solves the relation 
\begin{equation}
-\left( f,A_{D}v\right) _{\mathcal{O}}=l_{\phi }(v)\text{, for every }v\in
D(A_{D})\text{,}  \label{var}
\end{equation}%
where $l_{\phi }(\cdot )$ as given by (\ref{e3}) is a well-defined element
of $\left[ D(A_{D})\right] ^{\prime }$, by Step 3. Accordingly, we can use (%
\ref{iso}) and (\ref{e4}) to express the solution $f$ of (\ref{D})
as 
\[
\begin{array}{l}
f=-A_{D}^{-1}l_{\phi }\in L^{2}(\mathcal{O}); \\[.2cm] 
\left\Vert f\right\Vert _{L^{2}(\mathcal{O})}\leq \left\Vert
A_{D}^{-1}\right\Vert _{\mathcal{L}(\left[ D(A_{D})\right] ^{\prime },L^{2}(%
\mathcal{O}))}\left\Vert l_{\phi }\right\Vert _{\left[ D(A_{D})\right]
^{\prime }}\leq C_{2}\left\Vert \phi \right\Vert _{H^{-\frac{1}{2}}(\Omega
)}.%
\end{array}%
\]

\end{proof}

\bigskip

\subsubsection{The Argument for Dissipativity}\label{disscalc}

\label{diss} Considering the inner-product for the state space $\mathcal{H}$
given above in \eqref{innerproduct}, for any $\mathbf{y}=[p_{0},\mathbf{u}%
_{0},w_{0},w_{1}]\in D(\mathcal{A})$ we have 

\begin{eqnarray}
\left( \mathcal{A}\left[ 
\begin{array}{c}
p_{0} \\ 
\mathbf{u}_{0} \\ 
w_{0} \\ 
w_{1}%
\end{array}%
\right] ,\left[ 
\begin{array}{c}
p_{0} \\ 
\mathbf{u}_{0} \\ 
w_{0} \\ 
w_{1}%
\end{array}%
\right] \right) _{\mathcal{O}} &=&-\left( \mathbf{U}\cdot \nabla
p_{0},p_{0}\right) _{\mathcal{O}}-\left( \text{div}(\mathbf{u}%
_{0}),p_{0}\right) _{\mathcal{O}}  \nonumber \\
&&+\left( \text{div}\sigma (\mathbf{u}_{0})-\nabla p_{0}-\eta \mathbf{u}_{0}-%
\mathbf{U}\cdot \nabla \mathbf{u}_{0},\mathbf{u}_{0}-\alpha D(\mathbf{g}%
\cdot \nabla w_{0})\mathbf{e}_{3}\right) _{\mathcal{O}}  \nonumber \\
&&-\alpha \left( D(\mathbf{g}\cdot \nabla w_{1})\mathbf{e}_{3},\mathbf{u}%
_{0}-\alpha D(\mathbf{g}\cdot \nabla w_{0})\mathbf{e}_{3}\right) _{\mathcal{O%
}}  \nonumber \\
&&+\left( \Delta w_{1},\Delta w_{2}\right) _{\Omega }-\left( \Delta
^{2}w_{0},w_{1}+\mathbf{h}_{\alpha }\cdot \nabla w_{0}\right) _{\Omega } 
\nonumber \\
&&-\left( \left[ 2\nu \partial _{x_{3}}(\mathbf{u}_{0})_{3}+\lambda \text{div%
}(\mathbf{u}_{0})\right] _{\Omega }-\left. p_{0}\right\vert _{\Omega },w_{1}+%
\mathbf{h}_{\alpha }\cdot \nabla w_{0}\right) _{\Omega }  \nonumber \\
&&+\left( \mathbf{h}_{\alpha }\cdot \nabla w_{1},w_{1}+\mathbf{h}_{\alpha
}\cdot \nabla w_{0}\right) _{\Omega }  \label{diss_1}
\end{eqnarray}

Moreover, via Green's Theorem, as well as the assumption that $\mathbf{U}\in 
\mathbf{V}_{0}$ (as defined in (\ref{V_0})), we obtain 
\begin{equation}
2\text{Re}(\mathbf{U}\cdot \nabla p_{0},p_{0})_{\mathcal{O}}=-\int_{\mathcal{%
O}}\text{div}(\mathbf{U})\left\vert p_{0}\right\vert ^{2}d\mathcal{O};
\label{diss_2}
\end{equation}%
\begin{equation}
2\text{Re}(\mathbf{U}\cdot \nabla \mathbf{u}_{0},\mathbf{u}_{0})_{\mathcal{O}%
}=-\int_{\mathcal{O}}\text{div}(\mathbf{U})\left\vert \mathbf{u}%
_{0}\right\vert ^{2}d\mathcal{O}.  \label{diss_3}
\end{equation}

\noindent Applying Green's Theorem to right hand side of (\ref{diss_1}), with again $\mathbf{h}_{\alpha} = \mathbf U\big|_{\Omega}-\alpha \mathbf g$, we subsequently have 
\begin{align}
\left( {\mathcal{A}}%
\begin{bmatrix}
p_{0} \\ 
\mathbf{u}_{0} \\ 
w_{0} \\ 
w_{1}%
\end{bmatrix}%
,%
\begin{bmatrix}
p_{0} \\ 
\mathbf{u}_{0} \\ 
w_{0} \\ 
w_{1}%
\end{bmatrix}%
\right) _{\mathcal{H}}=& ~
 -(\sigma (\mathbf{u}_{0}),\epsilon (\mathbf{u}_{0}))_{\mathcal{O}}-\eta ||%
\mathbf{u}_{0}||_{\mathbf{L}^{2}(\mathcal{O})}^{2}   +\left\langle \sigma (\mathbf{u}_{0})\mathbf{n}-p_{0}%
\mathbf{n},\mathbf{u}_{0}\right\rangle _{\partial \mathcal{O}} \notag \\
&-(\left[ 2\nu \partial _{x_{3}}(\mathbf{u}_{0})_{3}+\lambda \text{div}(%
\mathbf{u}_{0})\right] _{\Omega }-\left. p_{0}\right\vert _{\Omega
},w_{1}+\mathbf U\cdot \nabla w_0)_{\Omega } \nonumber \\
& +\frac{1}{2}\int_{\mathcal O}\text{div}(\mathbf U)[|p_0|^2+|\bu_0|^2]d\mathcal O \notag \\
&-(\Delta^2w_0,\mathbf{ h}_{\alpha} \cdot \nabla w_0)_{\Omega}+(\mathbf{h}_{\alpha} \cdot \nabla w_1,w_1+\mathbf{h}_{\alpha} \cdot \nabla w_0)_{\Omega}\nonumber\\
&-\big(\text{div}~\sigma(\bu_0)-\eta\bu_0-\mathbf U\cdot \nabla \bu_0-\nabla p_0~,~ \alpha D(\mathbf g\cdot \nabla w_0) \mathbf{e}_3 \big  )_{\mathcal O} \nonumber \\
&+\big([2\nu\partial_{x_3}(\mathbf u_0)_3+\lambda \text{div}(\bu_0)-p_0]_{\Omega}~,~\alpha \mathbf g\cdot \nabla w_0\big)_{\Omega} \nonumber \\
& -\alpha\big(D (\mathbf g\cdot \nabla w_1) \mathbf{e}_3 ~,~\bu_0-\alpha D (\mathbf g\cdot \nabla w_0) \mathbf{e}_3 \big)_{\mathcal O}\nonumber \\
&+2i\text{Im}(p_0,\text{div}(\bu_0))_{\mathcal O}-i\text{Im}(\mathbf U\cdot \nabla p_0,p_0)_{\mathcal O}\nonumber \\
&~-i\text{Im}(\mathbf U\cdot \nabla \bu_0,\bu_0)_{\mathcal O}+2i\text{Im}(\Delta w_1,\Delta w_0)_{\Omega}.
\label{dissi_1.2}
\end{align}%
Invoking the boundary conditions (A.iv) and (A.v), in the definition of
the domain $D(\mathcal{A})$, there is then a cancellation
 of  terms at the boundary,
yielding
\begin{align}
\left( {\mathcal{A}}%
\begin{bmatrix}
p_{0} \\ 
\mathbf{u}_{0} \\ 
w_{0} \\ 
w_{1}%
\end{bmatrix}%
,%
\begin{bmatrix}
p_{0} \\ 
\mathbf{u}_{0} \\ 
w_{0} \\ 
w_{1}%
\end{bmatrix}%
\right) _{\mathcal{H}}=& ~ -(\sigma (\mathbf{u}_{0}),\epsilon (\mathbf{u}_{0}))_{\mathcal{O}}-\eta ||%
\mathbf{u}_{0}||_{\mathbf{L}^{2}(\mathcal{O})}^{2}   \notag \\
& +\frac{1}{2}\int_{\mathcal O}\text{div}(\mathbf U)[|p_0|^2+|\bu_0|^2]d\mathcal O \nonumber\\
&+I_1 + I_2 +I_3+I_4,
 \label{diss_1.4}
\end{align}
with terms grouped as
\begin{align}
I_1 = & ~-(\Delta^2w_0,\mathbf {h}_{\alpha} \cdot \nabla w_0)_{\Omega} \label{I1} \\
I_2=&~-\big(\text{div}~\sigma(\bu_0)-\eta\bu_0-\mathbf U\cdot \nabla \bu_0-\nabla p_0~,~\alpha D(\mathbf g\cdot \nabla w_0) \mathbf{e}_3 \big)_{\mathcal O}  \nonumber \\
&+\big([2\nu\partial_{x_3}(\mathbf u_0)_3+\lambda \text{div}(\bu_0)-p_0]_{\Omega}~,~\alpha \mathbf g\cdot \nabla w_0\big)_{\Omega} \label{I2} \\
I_3 = &~ -\alpha\big(D(\mathbf g\cdot \nabla w_1) \mathbf{e}_3 ~,~\bu_0-\alpha D(\mathbf g\cdot \nabla w_0) \mathbf{e}_3 \big)_{\mathcal O} \nonumber \\
&+(\mathbf {h}_{\alpha} \cdot \nabla w_1,w_1+\mathbf{h}_{\alpha} \cdot \nabla w_0)_{\Omega} \label{I3} \\
I_4 = & ~2i\text{Im}(p_0,\text{div}(\bu_0))_{\mathcal O}-i\text{Im}(\mathbf U\cdot \nabla p_0,p_0)_{\mathcal O}\nonumber \\
&~-i\text{Im}(\mathbf U\cdot \nabla \bu_0,\bu_0)_{\mathcal O}+2i\text{Im}(\Delta w_1,\Delta w_0)_{\Omega} \label{I4}
\end{align}
In order to have that $\widehat{\mathcal A}$ is dissipative in $\mathcal H$, we will control each term $ I_{i} $ $(i=1,...,4)$ by $||\mathbf y||^2_{\mathcal H}$. For the $I_1$ term, we follow the standard calculations typically used for the so called {\em flux multiplier} in the context of boundary control for the wave equation:
\begin{align}
-(\Delta^2w_0,\mathbf {h}_{\alpha} \cdot \nabla w_0)_{\Omega}=&~(\nabla \Delta w_0,\nabla(\mathbf {h}_{\alpha} \cdot \nabla w_0))_{\Omega}+0   \nonumber \\
=&~-(\Delta w_0,\Delta(\mathbf {h}_{\alpha} \cdot \nabla w_0))_{\Omega}+\int_{\partial \Omega}(\mathbf {h}_{\alpha} \cdot \mathbf \nu)|\Delta w_0|^2 d\partial \Omega, \label{mult_s}
\end{align}
where, in the first equality we have directly invoked the clamped plate boundary conditions, and in the second we have used the fact that $w_0=\partial_{\mathbf \nu}w_0=0$ on $\partial \Omega$ which yields that
$$ \frac{\partial}{\partial \mathbf{\nu}}(\mathbf {h}_{\alpha}  \cdot \nabla w_0)=(\mathbf {h}_{\alpha}  \cdot \mathbf{\nu})\frac{\partial^2 w_0}{\partial \mathbf{\nu}}=(\mathbf {h}_{\alpha}  \cdot \mathbf{\nu})(\Delta w_0\big|_{\partial \Omega}).  $$
 
\noindent (See \cite{lagnese} or \cite[p.305]{redbook}.)
Using the commutator bracket $[\cdot,\cdot]$, we can rewrite \eqref{mult_s} as
\begin{align}
-(\Delta^2w_0,\mathbf {h}_{\alpha} \cdot \nabla w_0)_{\Omega}=&~-(\Delta w_0,[\Delta,\mathbf {h}_{\alpha} \cdot \nabla]w_0)_{\Omega}-(\Delta w_0,\mathbf {h}_{\alpha} \cdot \nabla(\Delta w_0))_{\Omega}+\int_{\partial \Omega}(\mathbf {h}_{\alpha} \cdot \mathbf \nu)|\Delta w_0|^2d\partial\Omega.
\end{align}
With Green's relations once more:
\begin{align}
-(\Delta^2w_0,\mathbf {h}_{\alpha} \cdot \nabla w_0)_{\Omega}=&~-(\Delta w_0,[\Delta,\mathbf {h}_{\alpha} \cdot \nabla]w_0)_{\Omega}-\frac{1}{2}\int_{\partial \Omega}(\mathbf {h}_{\alpha} \cdot \mathbf \nu)|\Delta w_0|^2d\partial\Omega\nonumber \\
&+\frac{1}{2}\int_{\Omega}\big[\text{div}(\mathbf {h}_{\alpha} )\big]|\Delta w_0|^2d\Omega-i\text{Im}(\Delta w_0,\mathbf {h}_{\alpha} \cdot \nabla (\Delta w_0))_{\Omega} \nonumber \\
&+\int_{\partial \Omega} (\mathbf {h}_{\alpha} \cdot \mathbf \nu)|\Delta w_0|^2d\partial\Omega.
\end{align}
Thus, the final identity for $I_1$ is
\begin{align}
I_1=&~-(\Delta w_0,[\Delta,\mathbf {h}_{\alpha} \cdot \nabla]w_0)_{\Omega}+\frac{1}{2}\int_{\partial \Omega}(\mathbf {h}_{\alpha} \cdot \mathbf \nu)|\Delta w_0|^2d\partial\Omega\nonumber \\
&+\frac{1}{2}\int_{\Omega}\big[\text{div}(\mathbf {h}_{\alpha} )\big]|\Delta w_0|^2d\Omega-i\text{Im}(\Delta w_0,\mathbf {h}_{\alpha} \cdot \nabla (\Delta w_0)).
\end{align}
Now, recall that $\mathbf {h}_{\alpha}  =\mathbf U\big|_{\Omega}-\alpha \mathbf g$, where $\mathbf g$ is a $C^2(\Omega)$ -extension of exterior unit normal $\mathbf \nu(\xb)$, and $\alpha$ a parameter. If now, $\alpha$ satisfies
\begin{equation}\label{alphabound}
\alpha \ge \max_{\xb \in \partial \Omega} \mathbf U(\xb)\cdot \mathbf \nu(\xb), \end{equation}
we will have 
\begin{align}
\text{Re}~ I_1 =&~ \dfrac{1}{2}\int_{\partial \Omega} (\mathbf U\cdot \mathbf \nu-\alpha)|\Delta w_0|^2d\partial\Omega +\dfrac{1}{2}\int_{\Omega}\text{div}(\mathbf {h}_{\alpha} )|\Delta w_0|^2d\Omega - \text{Re}(\Delta w_0, [\Delta,\mathbf {h}_{\alpha} \cdot \nabla]w_0)_{\Omega}\nonumber \\
\le&~ 0 + \Big|\dfrac{1}{2}\int_{\Omega}\text{div}(\mathbf {h}_{\alpha} )|\Delta w_0|^2d\Omega - \text{Re}(\Delta w_0, [\Delta,\mathbf {h}_{\alpha} \cdot \nabla]w_0)_{\Omega} \Big|.  \label{I_1}
\end{align}
Since we can explicitly compute the commutator
\begin{align}
[\Delta, {\mathbf {h}_{\alpha} }\cdot \nabla]w_0 = & (\Delta h_1)(\partial_{x_1}w_0)+2(\partial_{x_1}h_1)(\partial_{x_1}^2w_0)\nonumber +2(\partial_{x_2} h_2)(\partial_{x_2}^2w_0)+(\Delta h_2)(\partial_{x_2}w_0) \\ 
&+2\text{div}({\mathbf {h}_{\alpha} })(\partial_{x_1}\partial_{x_2} w_0)
\end{align}
(where $\mathbf {h}_{\alpha} =[h_1,h_2]$), it is clear that 
\begin{equation}\label{commest}
\big|\big| [\Delta,{\mathbf {h}_{\alpha}}\cdot \nabla]w_0 \big|\big|_{L^2(\Omega)} \le C(\mathbf {h}_{\alpha}) ||\Delta w_0||_{L^2(\Omega)}.
\end{equation}
Hence, we get from \eqref{I_1},
\begin{equation}\label{I1bound} \text{Re}~I_1 \le C_1(\mathbf {h}_{\alpha})||\Delta w_0||^2_{\Omega}, \end{equation} for $\alpha$ sufficiently large. 

\medskip
 For the term $I_{2}$, as given in (\ref{I2}), we use Green's formula and the
fact that normal vector $\left. \mathbf{n}\right\vert _{\Omega}=[0,0,1]$, so as to have 
\begin{eqnarray}
I_{2} &=&\alpha \left( \sigma (\mathbf{u}_{0}),\epsilon (D(\mathbf{g}\cdot
\nabla w_{0})\mathbf{e}_{3})\right) _{\mathcal{O}}-\alpha \left( p_{0},\text{%
div}[D(\mathbf{g}\cdot \nabla w_{0})\mathbf{e}_{3}]\right) _{\mathcal{O}} 
\nonumber \\
&&+\alpha \left( \mathbf{U}\cdot \nabla \mathbf{u}_{0},D(\mathbf{g}\cdot
\nabla w_{0})\mathbf{e}_{3}\right) _{\mathcal{O}}+\eta \left( \mathbf{u}%
_{0},D(\mathbf{g}\cdot \nabla w_{0})\mathbf{e}_{3}\right) _{\mathcal{O}}.
\label{I2_2}
\end{eqnarray}

\smallskip
\noindent
Estimating right hand side by \eqref{D}--\eqref{D1}, Korn's inequality, and Young's inequality $\left\vert ab\right\vert \leq \delta a^{2}+C_{\delta }b^{2}$, we then have

\begin{equation}
\text{Re }I_{2}\leq \delta \left[ \left( \sigma (\mathbf{u}_{0}),\epsilon (%
\mathbf{u}_{0}))\right) _{\mathcal{O}}+\eta \left\Vert \mathbf{u}%
_{0}\right\Vert _{\mathcal{O}}^{2}\right] +C(\alpha ,\delta )\left[
\left\Vert p_{0}\right\Vert _{\mathcal{O}}^{2}+\left\Vert \Delta
w_{0}\right\Vert _{\mathcal{O}}^{2}\right] \text{, \ }\forall \text{ }\delta
>0  \label{I2bound}
\end{equation}

\medskip

Now we proceed with $I_3$, as given in \eqref{I3}. We firstly note that by condition (A.v) in the definition of $D(\mathcal A)$,
 
$$w_1=(\bu_0)_3-\mathbf U\cdot \nabla w_0~~~\text{ on }~~\Omega,$$ 
since $\mathbf n = [ 0,0,1]$ on $\Omega$. Thus we have the identity
$$\mathbf g\cdot \nabla w_1 = \mathbf g\cdot \nabla [(\bu_0)_3]-\mathbf g \cdot \nabla [\mathbf U\cdot \nabla w_0];$$
whence we obtain the estimate
\begin{align}
||\mathbf g\cdot \nabla w_1||_{H^{-1/2}(\Omega)} \le&~ C\Big[||(\bu_0)_3||_{H^{1/2}(\Omega)}+||\mathbf U \cdot \nabla w_0||_{H^1(\Omega)} \Big] \nonumber \\
\le & ~ C(\mathbf U)\Big[||\bu_0||_{\mathbf H^1(\mathcal O)}+||\Delta w_0||_{\Omega}\Big].
\end{align}
Using this estimate and the boundedness of the Dirichlet map $D$ in Lemma \ref{mark} for $H^{-1/2}(\Omega)$ boundary data, as well as \eqref{D}--\eqref{D1} for $H^{1/2 + \epsilon}_0(\Omega)$ boundary data, we have 
\begin{align}\label{using1}
\Big| \alpha\big(D(\mathbf g\cdot \nabla w_1) \mathbf{e}_3~,~\bu_0-\alpha D(\mathbf g\cdot \nabla w_0) \mathbf{e}_3 \big)_{\mathcal O} \Big| \le & ~C\Big[||\bu_0||_{\mathbf H^1(\mathcal O)}+||\Delta w_0||_{\Omega}\Big]||\bu_0||_{\mathcal O} \nonumber \\
&+C\Big[||\bu_0||_{\mathbf H^1(\mathcal O)}+||\Delta w_0||_{\Omega}\Big]||\Delta w_0||_{\Omega} \nonumber \\
\le & ~  \delta \Big[(\sigma(\bu_0),\epsilon(\bu_0))_{\mathcal O}+\eta ||\bu_0||^2\Big] \nonumber \\
&+C(\delta,\alpha)\Big[||\bu_0||^2_{\mathcal O}+||\Delta w_0||^2_{\Omega}\Big],
\end{align}
where we have again used Korn's and Young's inequalities. For the second term in $I_3$ in \eqref{I3}, owing to $w_1\big|_{\partial \Omega}=0$, we get
\begin{equation}
(\mathbf {h}_{\alpha} \cdot \nabla w_1,w_1)_{\Omega}=-\frac{1}{2}\int_{\Omega}\text{div}(\mathbf {h}_{\alpha})|w_1|^2d\Omega+i\text{Im}(\mathbf {h}_{\alpha}\cdot \nabla w_1,w_1)_{\Omega}.
\end{equation}
Additionally, since we have
\begin{equation}
(\mathbf {h}_{\alpha}\cdot \nabla w_1, \mathbf {h}_{\alpha}\cdot \nabla w_0)_{\Omega} = -(w_1, \text{div}(\mathbf {h}_{\alpha})\mathbf {h}_{\alpha}\cdot \nabla w_0+\mathbf {h}_{\alpha}\cdot \nabla [\mathbf {h}_{\alpha}\cdot \nabla w_0])_{\Omega},
\end{equation} then
\begin{align}\label{using2}
\text{Re}~(\mathbf {h}_{\alpha}\cdot \nabla w_1~,~w_1+\mathbf {h}_{\alpha} \cdot \nabla w_0)_{\Omega}  \le & ~ \Big|\frac{1}{2}\int_{\Omega} \text{div}(\mathbf {h}_{\alpha})|w_1|^2d\Omega +(w_1,\text{div}(\mathbf {h}_{\alpha})\mathbf {h}_{\alpha}\cdot \nabla w_0+\mathbf {h}_{\alpha}\cdot \nabla [\mathbf {h}_{\alpha}\cdot \nabla w_0])_{\Omega} \Big| \nonumber \\
\le & ~ C(\alpha)\big[||w_1||^2_{\Omega}+||\Delta w_0||^2_{\Omega}\big].
\end{align}
Using \eqref{using1} and \eqref{using2} in \eqref{I3}, we can estimate $I_3$ as
\begin{align}\label{I3bound}
\text{Re}~I_3 \le & ~ \delta\Big[(\sigma(\bu_0),\epsilon(\bu_0))_{\mathcal O}+\eta ||\bu_0||^2_{\mathcal O}\Big] + C(\alpha,\delta)\big[||\bu_0||_{\mathcal O}^2+||w_1||^2_{\Omega}+||\Delta w_0||^2_{\Omega}\big].
\end{align}
Recalling the definition of $\widehat{\mathcal A}$:
\begin{equation}
\widehat{\mathcal{A}}=\mathcal{A}-
\begin{bmatrix}
 \dfrac{\text{div}(\mathbf{U})}{2}+\varepsilon I & 0 & 0 & 0 \\ 
0 & \dfrac{\text{div}(\mathbf{U})}{2} +\varepsilon I & 0 & 0 \\ 
0 & 0 & \varepsilon I & 0 \\ 
0 & 0 & 0 & \varepsilon I%
\end{bmatrix},
\end{equation} 
then from \eqref{diss_1.4}, we have that
 for $\mathbf y=[p_0,\bu_0,w_0,w_1] \in D(\mathcal A)$,
\begin{align}
\text{Re}\left( \widehat{\mathcal{A}}%
\begin{bmatrix}
p_{0} \\ 
\mathbf{u}_{0} \\ 
w_{0} \\ 
w_{1}%
\end{bmatrix}%
,%
\begin{bmatrix}
p_{0} \\ 
\mathbf{u}_{0} \\ 
w_{0} \\ 
w_{1}%
\end{bmatrix}%
\right) _{\mathcal{H}}=& ~ -(\sigma (\mathbf{u}_{0}),\epsilon (\mathbf{u}_{0}))_{\mathcal{O}}-\eta ||%
\mathbf{u}_{0}||_{\mathcal{O}}^{2}  -\varepsilon||\mathbf y||^2_{\mathcal H}\nonumber\\
&+\text{Re}~I_1+\text{Re}~I_2+\text{Re}~I_3,
\end{align}
where $ \text{Re}~I_4=0  $. With the bounds on the real parts of $I_1, I_2, I_3$ given respectively in \eqref{I1bound} (with $\alpha$ sufficiently large), \eqref{I2bound}, and \eqref{I3bound}, we have
\begin{align}
\text{Re}~(\widehat{\mathcal A}\mathbf y, \mathbf y)_{\mathcal H} \le&~ \Big[-\varepsilon+C(\delta)\Big]||\mathbf y||_{\mathcal H}^2 \nonumber \\ & +(2\delta-1)\big[(\sigma(\bu_0),\epsilon(\bu_0))_{\mathcal O}+\eta||\bu_0||^2_{\mathcal O}\big],
\end{align} 
for all $\delta >0$, where we have suppressed dependence on $\mathbf U, \alpha, \Omega, \mathbf g$ in the constant $C$. Choosing $ \delta\le 1/2$,  then choosing the perturbation parameter $\varepsilon $ sufficiently large, we see that 
$$\text{Re}~(\widehat{\mathcal A}\mathbf y, \mathbf y)_{\mathcal H}  \le 0.$$
Thus $\widehat{\mathcal A}$ is dissipative on $\mathcal H$ with respect to the modified inner product  in \eqref{innerproduct}.

\subsection{Maximality}
\label{max} In this section we show the maximality property of the operator $%
\widehat{\mathcal{A}}$ on the space $\mathcal{H}$. To this end, we will need
to establish the associated 
\emph{range condition}, at least for parameter $\xi >0~~$%
sufficiently large. Namely, we must show 
\begin{equation}
Range(\xi I-\widehat{\mathcal{A}})=\mathcal{H},~~\text{for some}~~\xi >0.
\label{range_0}
\end{equation}%
This necessity is equivalent to finding $[p_0,{\mathbf{v}
_0},w_{1},w_{2}]\in D(%
\mathcal{A})$ which satisfies, for given $[p^{\ast },{\mathbf{v}}^{\ast
},w_{1}^{\ast },w_{2}^{\ast }]\in \mathcal{H}$, the abstract equation 
\begin{equation}
(\xi I-\widehat{\mathcal{A}})%
\begin{bmatrix}
p_0 \\ 
{\mathbf{v}_0} \\ 
w_{1} \\ 
w_{2}%
\end{bmatrix}%
=%
\begin{bmatrix}
p^{\ast } \\ 
{\mathbf{v}}^{\ast } \\ 
w_{1}^{\ast } \\ 
w_{2}^{\ast }%
\end{bmatrix}%
.  \label{range}
\end{equation}

Given the definition of $\mathcal{A}$ in (\ref{AAA}), and of $\widehat{\mathcal{A}}$ in (\eqref{pertA}),
solving the abstract equation (\ref{range}) is equivalent to proving that
the following system of equations, with given data $[p^{\ast },{\mathbf{v}}%
^{\ast },w_{1}^{\ast },w_{2}^{\ast }]\in \mathcal{H}$, has a (unique)
solution $[p_0,{\mathbf{v}_0},w_{1},w_{2}]\in D(\mathcal{A})$:%
\begin{align}
& \left\{ 
\begin{array}{l}
(\varepsilon+\xi) p_0+\mathbf{U}\cdot \nabla p_0+\frac{1}{2}\text{div}(\mathbf{U})p_0+\text{div}(%
{\mathbf{v}_0})=\text{ }p^{\ast }~~\text{ in }~\mathcal{O} \\ 
(\varepsilon+\xi) {\mathbf{v}_0}+\mathbf{U}\cdot \nabla {\mathbf{v}_0}+\frac{1}{2}\text{div}(%
\mathbf{U}){\mathbf{v}_0}-\text{div}~\sigma ({\mathbf{v}_0})+\eta {\mathbf{v}_0}%
+\nabla p_0=~{\mathbf{v}}^{\ast }~~\text{ in }~\mathcal{O} \\ 
(\sigma (\mathbf{v}_0)\mathbf{n}-p_0\mathbf{n})\cdot \boldsymbol{\tau }=0~\text{
on }~\partial \mathcal{O} \\ 
\mathbf{v}_0 \cdot \mathbf{n}=0~\text{ on }~S \\ 
\mathbf{v}_0 \cdot \mathbf{n}=w_{2}+\mathbf U\cdot \nabla w_1~\text{ on }~\Omega%
\end{array}%
\right.  \label{statice1} \\
& \left\{ 
\begin{array}{l}
(\varepsilon+\xi) w_{1}-w_{2}=~w_{1}^{\ast }~~\text{ in }~\Omega \\ 
(\varepsilon+\xi) w_{2}+\Delta ^{2}w_{1}+\left[ 2\nu \partial _{x_{3}}(\mathbf{v}_0)%
_{3}+\lambda \text{div}(\mathbf{v}_0)-p_0\right] _{\Omega }=w_{2}^{\ast }~\text{
in }~\Omega \\ 
w_{1}=\frac{\partial w_{1}}{\partial {\mathbf \nu} }=0~\text{ on }~\partial \Omega .%
\end{array}%
\right.  \label{staticsys3.5}
\end{align}%
We recall that {\em the parameter $\varepsilon >0$ is now fixed}, having been taken sufficiently large, in order that $\widehat{\mathcal A}$ be dissipative. 

The key ingredient of the following proof will be the
well-posedness result from \cite{AGW1} (itself based on \cite{dV}) which applies to (uncoupled)
equations of the type satisfied by the pressure variable. (Also see \cite{LaxPhil}.) We will then proceed to establish the range condition (\ref{range}), by
sequentially proving the existence of the pressure-fluid-structure
components $\left\{ p,\mathbf{v},w_{1},w_{2}\right\} $ which solve the
coupled system (\ref{statice1})--(\ref{staticsys3.5}). This work for
pressure-fluid-structure static well-posedness involves appropriate uses of
the Lax-Milgram Theorem. Now, let us give the following key lemma \cite{AGW1}, \cite{dV}:

\begin{lemma}
\label{staticwellp} Let $\mathbf{U}\in \mathbf{V}%
_{0}\cap \mathbf{H}^{3}(\mathcal{O})$ and consider the following $\xi $-parameterized PDE system on the fluid domain $%
\mathcal{O}$, with given forcing terms $\left\{ p^{\ast },{\mathbf{v}}^{\ast
}\right\} \in {L}^{2}(\mathcal{O})\times \lbrack \mathbf{V}_{0}]^{\prime }$%
\textbf{\ }and boundary data $g\in H_{0}^{1/2+\delta }(\Omega )$, where $%
\delta >0$:

\begin{align}
\xi p+\mathbf{U}\cdot \nabla p+\frac{1}{2}\text{div}(\mathbf{U})p+\text{div}(%
{\mathbf{v}})=& ~p^{\ast }~~\text{ in }~\mathcal{O}  \label{a1} \\
\xi {\mathbf{v}}+\mathbf{U}\cdot \nabla {\mathbf{v}}+\frac{1}{2}\text{div}(%
\mathbf{U}){\mathbf{v}}-\text{div}~\sigma ({\mathbf{v}})+\eta {\mathbf{v}}%
+\nabla p=& ~{\mathbf{v}}^{\ast }~~\text{ in }~\mathcal{O}  \label{a2} \\
\left( \sigma ({\mathbf{v}})\mathbf{n}-p\mathbf{n}\right) \cdot \boldsymbol{%
\tau }=& ~0~~\text{ on }~\partial \mathcal{O}  \label{a3} \\
{\mathbf{v}}\cdot \mathbf{n}=& ~0~~\text{ on }~S  \label{a4} \\
{\mathbf{v}}\cdot \mathbf{n}=& ~g~~\text{ on}~\Omega.  \label{a5}
\end{align}

\noindent Then for $\xi >0$ sufficiently large, there exists a unique solution $%
\left\{p,\mathbf{v}\right\} $ $\in {L}^{2}(\mathcal{O})\times \mathbf{H}^{1}(%
\mathcal{O})$ of (\ref{a1})--(\ref{a5}) such that:

(i) The fluid solution component ${\mathbf{v}}$ is of the form 
\begin{equation}
\mathbf{v}=\mathbf{u}+{\mathbf{v}}_{g}\text{,}  \label{crucial}
\end{equation}%
where $\mathbf{u}\in \mathbf{V}_{0}$, and ${\mathbf{v}}_{g}\in 
\mathbf{H}^{1}(\mathcal{O})$ satisfies%
\begin{equation}
{\mathbf{v}}_{g}\Big|_{\partial \mathcal{O}}=%
\begin{cases}
0 & ~\text{ on }~S \\ 
g\mathbf{n} & ~\text{ on }~\Omega .%
\end{cases}
\label{crucial2}
\end{equation}%
(ii) The trace term $\left[ \sigma ({\mathbf{v}})\mathbf{n}-p\mathbf{n}%
\right] _{\partial \mathcal{O}}\in \mathbf{H}^{-\frac{1}{2}}(\partial 
\mathcal{O})$, and moreover satisfies 
\begin{equation}
\left\langle \sigma ({\mathbf{v}})\mathbf{n}-p\mathbf{n,\tau }\right\rangle
_{\mathbf{H}^{-\frac{1}{2}}(\partial \mathcal{O})\times \mathbf{H}^{\frac{1}{%
2}}(\partial \mathcal{O})}=0\text{ \ for all }\mathbf{\tau }\in
TH^{1/2}(\partial \mathcal{O}),  \label{crucial3}
\end{equation}%
and so the boundary condition (\ref{a3}) is satisfied in the sense of
distributions; see (\ref{weak2}) of Remark \ref{weaker}.

(iii) The pressure and fluid solution components $(p,\mathbf{v})$ satisfies
the following estimates, for $\xi =\xi (\mathbf{U})$ large enough: 
\begin{eqnarray}
\left\Vert p\right\Vert _{L^2(\mathcal{O})} &\leq &\frac{C}{\xi }\left\Vert
[p^{\ast },{\mathbf{v}}^{\ast },g]\right\Vert _{\mathbf{L}^{2}(\mathcal{O}%
)\times \lbrack \mathbf{V}_{0}]^{\prime }\times H_{0}^{\frac{1}{2}+\delta
}(\Omega )};  \label{cdd_0} \\
\left\Vert {\mathbf{v}}\right\Vert _{\mathbf{H}^{1}(\mathcal{O})} &\leq
&C\left\Vert [p^{\ast },{\mathbf{v}}^{\ast },g]\right\Vert _{\mathbf{L}^{2}(%
\mathcal{O})\times \lbrack \mathbf{V}_{0}]^{\prime }\times H_{0}^{\frac{1}{2}%
+\delta }(\Omega )}.  \label{cdd}
\end{eqnarray}
\end{lemma}
Now, with Lemma \ref{staticwellp} in hand, we properly deal with the coupled
fluid-structure PDE system (\ref{statice1})--(\ref{staticsys3.5}). Our
solution here will be predicated on finding the structural variable $w_{1}$
which solves the $\Omega $-problem (\ref{staticsys3.5}). 

 By virtue of Lemma \ref{staticwellp}, for given data pressure and fluid data $(p^{\ast },{\mathbf{v}}^{\ast })\in L^{2}(\mathcal{O})\times \mathbf{L}%
^{2}(\mathcal{O})$ from (\ref{statice1}) and given 
boundary data $z\in H_{0}^{1/2 + \delta}(\Omega )$,  we
have that the following problem has a unique solution $\left\{ p=p(z;p^{\ast };%
\mathbf{v}^{\ast }), \mathbf{v}=\mathbf{v}(z;p^{\ast };\mathbf{v}^{\ast })\right\} $: 
\begin{equation}
\begin{array}{l}
(\varepsilon+\xi)p(z;p^{\ast };\mathbf{v}^{\ast })+\mathbf{U}\cdot \nabla p(z;p^{\ast };%
\mathbf{v}^{\ast })+\frac{1}{2}\text{div}(\mathbf{U})p(z;p^{\ast };\mathbf{v}%
^{\ast })+\text{div}({\mathbf{v}}(z;p^{\ast };\mathbf{v}^{\ast }))=~p^{\ast
}~~\text{ in }~\mathcal{O} \\ 
(\xi+\varepsilon) {\mathbf{v}}(z;p^{\ast };\mathbf{v}^{\ast })+\mathbf{U}\cdot \nabla {%
\mathbf{v}}(z;p^{\ast };\mathbf{v}^{\ast })+\frac{1}{2}\text{div}(\mathbf{U})%
{\mathbf{v}}(z;p^{\ast };\mathbf{v}^{\ast }) \\ 
\text{ \ \ \ \ \ \ \ \ \ \ \ \ \ \ \ \ \ \ \ }-\text{div}~\sigma ({\mathbf{v}%
}(z;p^{\ast };\mathbf{v}^{\ast }))+\eta {\mathbf{v}}(z;p^{\ast };\mathbf{v}%
^{\ast })+\nabla p(z;p^{\ast };\mathbf{v}^{\ast })=~{\mathbf{v}}^{\ast }~~%
\text{ in }~\mathcal{O} \\ 
\left( \sigma ({\mathbf{v}}(z;p^{\ast };\mathbf{v}^{\ast }))\mathbf{n}%
-p(z;p^{\ast };\mathbf{v}^{\ast })\mathbf{n}\right) \cdot \boldsymbol{\tau }%
=~0~~\text{ on }~\partial \mathcal{O} \\ 
{\mathbf{v}}(z;p^{\ast };\mathbf{v}^{\ast })\cdot \mathbf{n}=~0~~\text{ on }%
~S \\ 
{\mathbf{v}}(z;p^{\ast };\mathbf{v}^{\ast })\cdot \mathbf{n}=~z~~\text{ on}%
~\Omega .%
\end{array}
\label{staticsys1}
\end{equation}

\noindent We decompose the solution of the boundary value problm (\ref{staticsys1}) into two parts:%
\begin{eqnarray}
\mathbf{v}(z;p^{\ast };\mathbf{v}^{\ast }) &=&\mathbf{v}(z)+\mathbf{v}(p^{\ast };%
\mathbf{v}^{\ast });  \label{d1} \\
p(z;p^{\ast };\mathbf{v}^{\ast }) &=&p(z)+p(p^{\ast };\mathbf{v}^{\ast }),
\label{d2}
\end{eqnarray}
where $(p(z),{\mathbf{v}}(z))\in  L^{2}(%
\mathcal{O})\times \mathbf{H}^{1}(\mathcal{O})$ is the solution of the problem 
\begin{align}
(\varepsilon+\xi)p(z)+\mathbf{U}\cdot \nabla p(z)+\frac{1}{2}\text{div}(\mathbf{U})p(z)+%
\text{div}({\mathbf{v}}(z))=& 0~~\text{ in }~\mathcal{O}  \label{huh} \\
(\varepsilon+\xi) {\mathbf{v}(z)}+\mathbf{U}\cdot \nabla {\mathbf{v}}(z)+\frac{1}{2}\text{%
div}(\mathbf{U}){\mathbf{v}}(z)-\text{div}~\sigma ({\mathbf{v}}(z))+\eta {%
\mathbf{v}}(z)+\nabla p(z)=& ~0~~\text{ in }~\mathcal{O} \\
\left( \sigma ({\mathbf{v}}(z))\mathbf{n}-p(z)\mathbf{n}\right) \cdot 
\boldsymbol{\tau }=& ~0~~\text{ on }~\partial \mathcal{O} \\
{\mathbf{v}}(z)\cdot \mathbf{n}=& ~0~~\text{ on }~S \\
{\mathbf{v}}(z)\cdot \mathbf{n}=& ~z~~\text{ on}~\Omega ;  \label{givensys2}
\end{align}%
and $\big(p(p^{\ast },{\mathbf{v}}^{\ast }),{\mathbf{v}}(p^{\ast };{\mathbf{v%
}}^{\ast })\big)\equiv (\overline{p},\overline{{\mathbf{v}}})\in L^{2}(%
\mathcal{O})\times \mathbf{V}_{0}$ ~~is the solution of the problem 
\begin{align}
(\varepsilon+\xi) \overline{p}+\mathbf{U}\cdot \nabla \overline{p}+\frac{1}{2}\text{div}(%
\mathbf{U})\overline{p}+\text{div}\overline{{\mathbf{v}}}=& ~p^{\ast }~\text{
in }~~\mathcal{O}  \label{insert} \\
(\varepsilon+\xi) \overline{{\mathbf{v}}}+\mathbf{U}\cdot \nabla \overline{{\mathbf{v}}}+%
\frac{1}{2}\text{div}(\mathbf{U})\overline{{\mathbf{v}}}-\text{div}~\sigma (%
\overline{{\mathbf{v}}})+\eta \overline{{\mathbf{v}}}+\nabla \overline{p}=& ~%
{\mathbf{v}}^{\ast }~~\text{ in }~~\mathcal{O} \\
\left( \sigma ({\overline{{\mathbf{v}}}})\mathbf{n}-\overline{p}\mathbf{n}\right) \cdot 
\boldsymbol{\tau }=& ~0~~\text{ on }~\partial \mathcal{O} \\
\overline{{\mathbf{v}}}\cdot \mathbf{n}=& ~0~~\text{ on }~~S \\
\overline{{\mathbf{v}}}\cdot \mathbf{n}=& ~0~~\text{ on }~~\Omega .
\label{insert1}
\end{align}%
Now, with $(p(z),{\mathbf{v}}(z))$ and $ (\overline{p},\overline{{\mathbf{v}}}) $ in hand,  if we multiply the structural PDE component (\ref{staticsys3.5}) by given $z\in
H_{0}^{2}(\Omega )$, invoke the resolvent relation
$w_2=(\varepsilon+\xi)w_1-w_1^*$, integrate by parts, and utilize
the boundary conditions in the BVP (\ref{huh})--(\ref{givensys2}) and (\ref%
{insert})--(\ref{insert1}), we then
have
\begin{align}
(\varepsilon+\xi)^2( w_1,z)_{\Omega}+( \Delta w_1,\Delta z)_{\Omega} +\langle \sigma ({\mathbf{v}_0})\mathbf{n}-p_0\mathbf{n},{\mathbf{v}}%
(z)\rangle _{\partial \mathcal{O}} = ( w_2^*+[\varepsilon+\xi]w_1^*,z)_{\Omega}.
\end{align}%
(In obtaining this relation that fluid solution component $v(z)$ of the BVP \eqref{insert}--(\ref{insert1} has the decomposition \eqref{crucial}-\eqref{crucial2}, with therein $g =z$.) A subsequent use of Green's identity and substitution then yields
\begin{align}
( w_{2}^{\ast }+[\xi+\varepsilon] w_{1}^{\ast },z) _{\Omega }=& ~(\xi+\varepsilon) ^{2}( w_{1},z) _{\Omega }+( \Delta w_{1},\Delta
z) _{\Omega }+(\varepsilon+\xi) ({\mathbf{v}_0},{\mathbf{v}}(z))_{\mathcal{O}}+(%
\mathbf{U}\cdot \nabla {\mathbf{v}_0},\mathbf{v}(z))_{\mathcal{O}}  \notag \\
& +(\sigma ({\mathbf{v}_0}),\epsilon ({\mathbf{v}}(z)))_{\mathcal{O}}+\frac{1}{%
2}(\text{div}(\mathbf{U}){\mathbf{v}_0},{\mathbf{v}}(z))_{\mathcal{O}}  \notag
\\
& +\eta ({\mathbf{v}_0},{\mathbf{v}}(z))_{\mathcal{O}}-(p_0,\text{div}({\mathbf{v%
}}(z)))_{\mathcal{O}}-({\mathbf{v}}^{\ast },{\mathbf{v}}(z))_{\mathcal{O}}.
\label{almost}
\end{align}%
Invoking
the respective solution maps for (\ref{huh})--(\ref{givensys2}) and (\ref%
{insert})--(\ref{insert1}), we may express the (prospective) solution
component $(p_0,{\mathbf{v}_0})$ of (\ref{statice1}) as 
\begin{align}
p_0& =p([\xi+\varepsilon] w_{1}+\mathbf U\cdot \nabla w_1-w_{1}^{\ast })+\overline{p} 
\label{r1} \\
 {\mathbf{v}_0}& ={\mathbf{v}}([\xi+\varepsilon] w_{1}+\mathbf U\cdot \nabla w_1-w_{1}^{\ast })~+\overline{{\mathbf{v}}}\label{r2}
\end{align}%
(cf. (\ref{d1})--(\ref{d2})). With (\ref{almost}) and (\ref{r1})--(\ref{r2})
in mind, we define an operator $B\in \mathcal{L}%
(H_{0}^{2}(\Omega ),H^{-2}(\Omega ))$ by:~~For $w$ and $z \in H_0^2(\Omega)$
\begin{align}
\left\langle B(w),z\right\rangle _{\big(H^{-2}(\Omega),H_{0}^{2}(\Omega )\big)}\equiv & ~(\varepsilon+\xi)^{2}(w,z) _{\Omega }+(
\Delta w,\Delta z) _{\Omega } +(\varepsilon+\xi)^2({\mathbf{v}}(w),{\mathbf{v}}(z))_{\mathcal{O}}  \nonumber \\ 
&+(\xi+\varepsilon)\big(\sigma ({\mathbf{v}}(w)),\epsilon ({\mathbf{v}}(z))\big)_{%
\mathcal{O}}+\eta (\xi+\varepsilon) \big({\mathbf{v}}(w),{\mathbf{v}}(z))_{\mathcal{O}%
}\nonumber \\ 
& +(\xi+\varepsilon)(\mathbf{U}%
\cdot \nabla {\mathbf{v}}(w),{\mathbf{v}}(z))_{\mathcal{O}}+\frac{\xi+\varepsilon }{2}%
\big(\text{div}(\mathbf{U}){\mathbf{v}}(w),{\mathbf{v}}(z)\big)_{\mathcal{O}}
\notag \\
& +(\xi+\varepsilon)\big(\mathbf v(\mathbf U\cdot \nabla w),\mathbf v(z)\big)_{\mathcal O}-\big(p(\mathbf U\cdot \nabla w),\text{div}(\mathbf v(z))\big)_{\mathcal O} \notag\\
& -(\xi+\varepsilon) (p(w),\text{div}({\mathbf{v}}(z)))_{\mathcal{O}}+B_0(w,z),\label{AA}
\end{align}%
where $(p(.),{\mathbf{v}}(.))$ solves (\ref{huh})--(\ref{givensys2}) for given 
$H_{0}^{2}(\Omega )$ boundary data, and $B_0: H_0^2(\Omega)\times H_0^2(\Omega) \to \mathbb C$ is given by
\begin{align}
B_0(w,z) = & ~\big( \mathbf U \cdot \nabla \mathbf v(\mathbf U\cdot \nabla w),\mathbf v(z)\big)_{\mathcal O} \notag\\
&+\frac{1}{2} \big(\text{div}(\mathbf U)\mathbf v(\mathbf U\cdot \nabla w),\mathbf v(z) \big)_{\mathcal O} + \eta \big( \mathbf v(\mathbf U\cdot \nabla w),\mathbf v(z)\big)_{\mathcal O}\notag\\
&+\big((\sigma(\mathbf v(\mathbf U\cdot \nabla w)),\epsilon(\mathbf v(z))\big)_{\mathcal O}. \notag
\end{align}
Then writing the relation (\ref{almost}) again and finding a solution $w_{1}\in
H_{0}^{2}(\Omega )$ for the structural PDE component (\ref{staticsys3.5}) is tantamount to finding solution $w_{1}\in H_{0}^{2}(\Omega )$
of the variational equation 
\begin{equation}\label{vareq}
\left\langle B(w_{1}),z\right\rangle _{\big(H^{-2}(\Omega),H_{0}^{2}(\Omega )\big)}=%
\mathcal{F}(z)\text{, for all }z\in H_{0}^{2}(\Omega )\text{,}  
\end{equation}%
where the functional $\mathcal{F}\in H^{-2}(\Omega )$ is given by 
\begin{align*}
\mathcal{F}(z)\equiv & ~( w_{2}^{\ast }+(\varepsilon+\xi) w_{1}^{\ast },z) _{\Omega } +\big({\mathbf{v}}^{\ast },{\mathbf{v}}(z)\big)_{%
\mathcal{O}}\\
& +(\xi +\varepsilon) ({\mathbf{v}}(w_{1}^{\ast }),{\mathbf{v}}(z))_{\mathcal{O}}-(\varepsilon+\xi) (%
\overline{{\mathbf{v}}},{\mathbf{v}}(z))_{\mathcal{O}} \\
& +(\mathbf{U}\cdot \nabla {\mathbf{v}}(w_{1}^{\ast }),{\mathbf{v}}(z))_{%
\mathcal{O}}-(\mathbf{U}\cdot \nabla \overline{{\mathbf{v}}},{\mathbf{v}}%
(z))_{\mathcal{O}} \\
& +\frac{1}{2}(\text{div}(\mathbf{U}){\mathbf{v}}(w_{1}^{\ast }),{\mathbf{v}}%
(z))_{\mathcal{O}}-\frac{1}{2}(\text{div}(\mathbf{U})\overline{{\mathbf{v}}},%
{\mathbf{v}}(z))_{\mathcal{O}} \\
& +\eta ({\mathbf{v}}(w_{1}^{\ast }),{\mathbf{v}}(z))_{\mathcal{O}}-\eta (%
\overline{{\mathbf{v}}},{\mathbf{v}}(z))_{\mathcal{O}} \\
& +(\sigma ({\mathbf{v}}(w_{1}^{\ast })),\epsilon ({\mathbf{v}}(z)))_{%
\mathcal{O}}-(\sigma (\overline{{\mathbf{v}}}),\epsilon ({\mathbf{v}}(z)))_{%
\mathcal{O}} \\
& -(p(w_{1}^{\ast }),\text{div}({\mathbf{v}}(z)))_{\mathcal{O}}+(\overline{p}%
,\text{div}({\mathbf{v}}(z)))_{\mathcal{O}}.
\end{align*}

\noindent We assert that the operator 
$B$ is $H_{0}^{2}(\Omega )$-elliptic, (so the relation (\ref{vareq}) can be solved by the Lax-Milgram Theorem) for $\xi =\xi (\mathbf{U})$ large
enough. To establish this fact, we start with the non-elliptic portion of $ B $ in (\ref{AA}). Let $w \in H_0^2(\Omega)$ be given. Then, by H\"older-Young inequalities we have 
\begin{align}\label{1}
(\varepsilon+\xi)\big(\mathbf U\cdot \nabla \mathbf v(w),\mathbf v(w)\big)_{\mathcal O} \ge&~ -\delta (\varepsilon+\xi)||\mathbf v(w)||^2_{\mathbf{H}^1(\mathcal O)}-C_{\delta}(\varepsilon+\xi)||\mathbf v(w)||^2_{\mathcal O}. 
\end{align}
In addition, via (\ref{cdd_0})--(\ref{cdd}) we have
\begin{align}
-(\varepsilon+\xi)\big(p(w),\text{div}(\mathbf v(w))\big)_{\mathcal O} \ge &~-(\varepsilon+\xi)\dfrac{C}{\xi+\varepsilon}||w||_{H_0^2(\Omega)}||\mathbf v(w)||_{\mathbf{H}^1(\mathcal O)} \nonumber \\
\ge & ~-C_{\delta}||\mathbf v(w)||^2_{\mathbf H^1(\mathcal O)}-\delta||\Delta w||^2_{\Omega}.  \label{2}
\end{align}
Moreover,
\begin{align}
(\xi+\varepsilon)\big(\mathbf v(\mathbf U\cdot \nabla w),\mathbf v(w)\big)_{\mathcal O} \ge &~-C_{\delta}||w||^2_{H^{3/2+\delta}_0(\Omega)}-\delta (\varepsilon+\xi)^2||\mathbf v(w)||^2_{\mathcal O} \nonumber \\
\ge & ~-C_{\delta}||w||^{1/2-\delta}_{\Omega}||\Delta w||^{3/2+\delta}_{\Omega}-\delta (\varepsilon+\xi)^2||\mathbf v(w)||^2_{\mathcal O} \nonumber\\
~\ge &~-C_{\delta}||w||^2_{\Omega}-\delta||\Delta w||^2_{\Omega}-\delta(\xi+\varepsilon)^2||\mathbf v(w)||^2_{\mathcal O}. \label{3}
\end{align}
And also, by (\ref{cdd_0}),
\begin{align}
-\big(p(\mathbf U\cdot \nabla w),\text{div}(\mathbf v(w))\big)_{\mathcal O} \ge &~-\dfrac{C}{(\xi+\epsilon)^2}||\Delta w||^2_{\Omega}-C||\mathbf v(w)||^2_{\mathbf H^1(\mathcal O)}. \label{4}
\end{align}
Standard interpolation and estimate (\ref{cdd}) give further that
\begin{align}
B_0(w,w) \ge & ~-C||w||_{H_0^1(\Omega)}||\Delta w||_{\Omega} \nonumber \\
\ge &~ -C||w||^{1/2}_{(\Omega)}||\Delta w||^{1/2}_{\Omega}||\Delta w||_{\Omega}  \nonumber \\
\ge &~-C_{\delta}||w||^2_{\Omega}-\delta||\Delta w||^2_{\Omega}. \label{5}
\end{align}
Additionally, it is clear that
\begin{align}
\dfrac{\xi + \varepsilon}{2}\big(\text{div}(\mathbf U)\mathbf v(w),\mathbf v(w)\big)_{\mathcal O} \ge&~-C_{\mathbf U}(\xi+\varepsilon)||\mathbf v(w)||^2_{\mathcal O}\label{6}.
\end{align}
Taking into account \eqref{1}--\eqref{6} in the relation \eqref{AA}, we have for $w\in H_0^2(\Omega)$ 
\begin{align}
\langle B(w),w\rangle_{(H^{-2}(\Omega),H_0^2(\Omega))} \ge &~\big[1- 3 \delta-C/(\xi+\epsilon)^2\big]||\Delta w||^2_{\Omega}\nonumber \\
&+ \big[(\xi+\varepsilon)^2-2C_{\delta}\big]||w||^2_{\Omega}\nonumber \\
&+\big[(\xi+\varepsilon)(1-\delta C_{\eta})-C_{\delta,\eta}]\Big\{(\sigma(\mathbf v(w),\epsilon(\mathbf v(w)))_{\mathcal O}+\eta||\mathbf v(w)||^2_{\mathcal O} \Big\} \nonumber \\
&+\big[(\xi+\varepsilon)^2(1-\delta)-(C_{\mathbf U}+C_{\delta})(\xi+\varepsilon)\big] ||\mathbf v(w)||^2_{\mathcal O} \nonumber \\
\ge & ~ C^*||\Delta w||^2_{\mathcal O}.
\end{align}
Here, $ C_{\eta} $ is the constant from Korn's inequality which has implicitly been used, and positive constant $C^*$ is independent of $\xi >0$, sufficiently large. 

Consequently, by Lax-Milgram Theorem there exists a unique solution $ w_1 $ to the variational
equation (\ref{vareq}), or what is the same, we can recover the solution component $%
w_{1}$ of the resolvent equations (\ref{statice1})--(\ref{staticsys3.5}). In
turn, we reconstruct the other solution variables of (\ref{statice1})--(\ref{staticsys3.5}) via
\begin{equation*}
\begin{array}{c}
{\mathbf{v}_0}=~{\mathbf{v}}((\xi+\varepsilon) w_{1}+\mathbf U\cdot \nabla w_1-w_{1}^{\ast })+\overline{{\mathbf{v}}}, \\ 
p_0 =~p((\xi+\varepsilon) w_{1}+\mathbf U\cdot \nabla w_1-w_{1}^{\ast })+\overline{p}\\
w_{2}=~(\xi+\varepsilon) w_{1}-w_{1}^{\ast },
\end{array}%
\end{equation*}%
where $(p(\cdot),\mathbf{v}(\cdot))$ refers to the solution to (\ref{huh}%
)--(\ref{givensys2}) with given boundary data, and $[\overline p,\overline{{\mathbf{v}}}]$ solves the
system (\ref{insert})--(\ref{insert1}).

Having recovered the solution $[p,\mathbf v, w_1, w_2]$, with the given data $[p^*,\mathbf v^*,w_1^*,w_2^*] \in \mathcal H$, we note that, a posteriori, the solution in fact resides in $\mathcal D(\mathcal A)$. Indeed, by Lax-Milgram, $w_1,w_2 \in H_0^2(\Omega)$, which immediately gives condition (A.iii) for $D(\mathcal A)$; with the solution $[p,\mathbf v, w_1, w_2] \in \mathcal H$, it is clear that having $\mathbf v \in \mathbf H^{1}(\mathcal O)$ and $p \in L^2(\mathcal O)$ gives that $\mathbf U\cdot \nabla p \in L^2(\mathcal O)$ which provides condition (A.i) for $D(\mathcal A)$. Condition (A.ii) then follows from the $\mathbf v$-equation in \eqref{statice1}, with (A.iv) and (A.v) coming from our use of Lemma \ref{staticwellp} to construct the solution.

This finally establishes the range condition in (\ref{range_0}) for $\xi >0$
sufficiently large. A subsequent application of Lumer-Philips Theorem yields
a contraction semigroup for the $\widehat{\mathcal{A}}:D(\mathcal{A})\subset 
\mathcal{H}\rightarrow \mathcal{H}$. As a consequence, the application of
Theorem 1.1 \cite[Chapter 3.1]{pazy}, p.76, gives the desired result for the
(unperturbed) compressible flow-structure generator $\mathcal{A}$.

\section{Acknowledgments}

The authors would like to sincerely thank Earl Dowell for his insight and expertise in discussing fluid flow modeling and coupling conditions \eqref{coupling}, as presented in Section \ref{modeling} \cite{dowellcorr}. 


The authors would like to thank the National Science
Foundation, and acknowledge their partial funding from NSF Grant
DMS-1616425 (G. Avalos and Pelin G. Geredeli). Pelin G. Geredeli also would like to thank the University of Nebraska-Lincoln for the Edith T. Hitz Fellowship.

\end{document}